\documentclass[12pt, reqno]{amsart}
\usepackage{pifont}
\usepackage{mathrsfs}
\usepackage{geometry}
\usepackage{titletoc}
\usepackage{stix2}
\usepackage{amscd}

\usepackage{amsmath}
\usepackage{amssymb} 
\usepackage{enumitem} 
\usepackage{mathtools} 
\usepackage[table]{xcolor} 
\usepackage[all]{xy} 
\usepackage{tikz} 
\usepackage{tikz-cd}
\usepackage{indentfirst} 
\usepackage{babel} 
\usepackage{setspace} 

\usepackage[colorlinks,linkcolor=red,anchorcolor=green,citecolor=blue]{hyperref} 
\hypersetup{linktocpage = true} 

\usepackage{rotating} 

\usepackage{ytableau} 
\usepackage{longtable} 
\newcolumntype{M}[1]{>{\centering\arraybackslash}m{#1}} 

\usepackage{fancyhdr}

\newcommand\cE{{\mathcal E}}

\newcommand\cL{{\mathcal L}}

\newcommand\cM{{\mathcal M}}

\newcommand\sF{{\mathscr F}}

\newcommand{\id}{{\rm id}}
\newcommand{\rk}{{\rm rk}}

\DeclareMathOperator*{\supp}{Supp}

\newcommand\bp{{\bar\partial}}

\usepackage{amsthm}
\theoremstyle{plain}

\newtheorem{theorem}{Theorem}[section]
\newtheorem{thm}{Theorem}[section]
\newtheorem{lemma}[thm]{Lemma}
\newtheorem{prop}[thm]{Proposition}
\newtheorem{cor}[thm]{Corollary}
\newtheorem{defn}[thm]{Definition}

\theoremstyle{definition}
\newtheorem{example}[thm]{Example}
\newtheorem{remark}[thm]{Remark}

\newcommand{\btheorem}{\begin{thm}}
    \newcommand{\etheorem}{\end{thm}}
\newcommand{\bproposition}{\begin{prop}}
    \newcommand{\eproposition}{\end{prop}}
\newcommand{\bdefinition}{\begin{defn}}
    \newcommand{\edefinition}{\end{defn}}
\newcommand{\bcorollary}{\begin{cor}}
    \newcommand{\ecorollary}{\end{cor}}
\newcommand{\bproof}{\begin{proof}}
    \newcommand{\eproof}{\end{proof}}
\newcommand{\bremark}{\begin{remark}}
    \newcommand{\eremark}{\end{remark}}
\newcommand{\eexample}{\end{example}}
\newcommand{\bexample}{\begin{example}}

\newcommand{\elemma}{\end{lemma}}
\newcommand{\blemma}{\begin{lemma}}

\newcommand{\sq}{\sqrt{-1}}

\newcommand{\p}{\partial}

\renewcommand{\bar}{\overline}
\newcommand{\eps}{\varepsilon}

\renewcommand{\phi}{\varphi}

\newcommand{\beq}{\begin{equation}}
\newcommand{\eeq}{\end{equation}}
\newcommand{\ee}{\end{eqnarray*}}
\newcommand{\be}{\begin{eqnarray*}}

\newcommand{\bd}{\begin{enumerate}}
    \newcommand{\ed}{\end{enumerate}}

\renewcommand{\tilde}{\widetilde}

\newcommand{\qtq}[1]{\quad\mbox{#1}\quad}
\renewcommand{\bp}{\bar{\partial}}
\newcommand{\Om}{\Omega}

\newcommand{\ts}{\otimes}

\renewcommand{\>}{\rightarrow}

\newcommand{\C}{{\mathbb C}}

\newcommand{\R}{{\mathbb R}}

\newcommand{\LL}{\left\langle}
\newcommand{\RL}{\right\rangle}

\newcommand{\Ric}{\mathrm{Ric}}

\renewcommand{\>}{\rightarrow}

\renewcommand{\p}{{\partial}}
\renewcommand{\bp}{{\bar{\partial}}}

\newcommand{\vone}{ \vskip 1\baselineskip}
\newcommand{\om}{\omega}

\renewcommand{\bar}{\overline}
\renewcommand{\tilde}{\widetilde}

\newcommand{\smo}{\sqrt{-1}}

\newcommand{\Herm}{\mathrm{Herm}}
\newcommand{\tr}{\mathrm{tr}}

\setlist[itemize]{leftmargin=*}
\setlist[enumerate]{leftmargin=*}

\numberwithin{equation}{section} 

\makeatletter

\usepackage{fancyhdr}
\title{The prescribed Hermitian-Yang-Mills flow I}

\author{Zhiyao Xiong}
\author{Xiaokui Yang}
\author{Shing-Tung Yau}

\address{Zhiyao Xiong, Department of Mathematics, Tsinghua University, Beijing, 100084, China}
\email{xiongzy22@mails.tsinghua.edu.cn}

\address{Xiaokui Yang, Department of Mathematics and Yau Mathematical Sciences Center, Tsinghua University, Beijing, 100084, China}
\email{xkyang@mail.tsinghua.edu.cn}

\address{Shing-Tung Yau, Yau Mathematical Sciences Center and  Qiuzhen College, Tsinghua University, Beijing, 100084, China}
\email{styau@mail.tsinghua.edu.cn}

\begin{document}

    \begin{abstract} In this paper, we introduce a broad class of flows, including the prescribed Hermitian-Yang-Mills flow:
    $$\frac{\p h}{\p t}=-\Lambda_{\omega_g}\left(\sq R^h\right)+P$$
   where $P\in\Gamma(M,E^*\otimes\bar{E}^*)$ is a prescribed Hermitian tensor associated with a holomorphic vector bundle $E$ over a  K\"ahler (or Hermitian) manifold $(M,\omega_g)$.  We establish the long-time convergence of the flow to a limiting metric $h_{\infty}$ and use it to solve the prescribed Hermitian-Yang-Mills tensor equation
    $$\Lambda_{\omega_g}\left(\sq R^{h_\infty}\right)=P,
    $$
    for a general class of prescribed Hermitian tensors $P$.
    The crucial uniform $C^0$-estimate of $\{h(t)\}$ along the flow is obtained via a parabolic comparison principle.
    \end{abstract}

    \maketitle {
        \setcounter{tocdepth}{1}

    {\small{    \begin{spacing}{1.1} \tableofcontents %
                \dottedcontents{section}[1.8cm]{}{3em}{5pt} %
\end{spacing} }} }

\vskip -1\baselineskip

    \section{Introduction}

    On a compact K\"ahler manifold $(M, \om_g)$, the Ricci form $\Ric(\om_g)$ represents the first Chern class $c_1(M)$ up to a factor of $2\pi$:
    \[
    [\Ric(\om_g)] = 2\pi c_1(M) \in H^{1,1}(M, \R).
    \]
    The celebrated solution to the Calabi conjecture by Yau \cite{Cal57,Yau78} establishes the converse--namely, the existence of metrics with prescribed Ricci curvature--forging a definitive link between geometry and topology:

    \setcounter{theorem}{0}
    \renewcommand{\thetheorem}{\Alph{theorem}}
        \begin{theorem}\label{thm:CY_intro}
        Let $(M, \om_g)$ be a compact K\"ahler manifold. For any closed real $(1,1)$-form $\Omega$ representing $2\pi c_1(M)$, there exists a unique K\"ahler metric $\omega \in [\om_g]$ such that
        \begin{equation}
        \Ric(\omega) = \Omega. \label{CY1}
        \end{equation}
    \end{theorem}

\noindent   Note  that, since $\omega$ is K\"ahler, the Calabi-Yau equation \eqref{CY1} is equivalent to
\beq \Lambda_{\om}\left(\sqrt{-1}\, R^{\omega}\right) =\Om.\label{CY}\eeq

    \noindent The Calabi-Yau theorem  has become a cornerstone of modern geometry, profoundly impacting differential geometry, algebraic geometry, and mathematical physics. Its significance extends beyond the K\"ahler setting to general Hermitian manifolds, where analogous existence results have been developed; see, e.g., \cite{FY08, TW10, STW17} and the extensive references therein.  \\

  Let's recall the classical Donaldson--Uhlenbeck--Yau theorem (\cite{Don85, UY86, Don87}) for stable vector bundles (see also the Hermitian analogue in \cite{LY86}):  \begin{theorem}  For  a stable holomorphic vector bundle $E$ over a compact K\"ahler manifold $(M,\omega_g)$,  there exists a unique Hermitian-Einstein metric $h$ on $E$ up to scaling, satisfying
    \beq \Lambda_{\omega_g} \left(\sq R^h\right)=\lambda_E\cdot h, \label{HE} \eeq
    where $R^h\in\Gamma(M,\Lambda^{1,1}T^*M\ts E^*\ts \bar E^*)$ is the Chern curvature tensor of $(E,h)$ and
    $$\lambda_E=\frac{2\pi n\int_M c_1(E)\wedge \omega_g^{n-1}}{\mathrm{rank}(E)\int_M\omega_g^n}.$$ \end{theorem}

\noindent This foundational result established a deep equivalence on compact K\"ahler manifolds: a holomorphic vector bundle admits a Hermitian--Einstein metric if and only if it is polystable in the sense of Mumford--Takemoto. This forged a critical bridge between differential geometry---the existence of special metrics solving a nonlinear PDE---and algebraic geometry, translating the purely algebraic notion of stability into an analytic criterion. The correspondence not only unified two seemingly disparate fields but also laid the essential groundwork for the development of complex geometry, gauge theory and representation theory.  For comprehensive technical details, we refer the reader to \cite{NS65, Hit87, Sim88, Sim92} and the extensive references cited therein.\\

In the series of papers \cite{WYY26+,FWYY26+,WYY26b+}, the second and third authors, together with M.-W. Wang and J.-X. Fang, established a vector bundle analogue of the Calabi--Yau theorem. This result provides a solution to the prescribed Hermitian--Yang--Mills tensor equation, extending the classical Hermitian--Einstein existence theory to the case of fully prescribed Hermitian-Yang-Mills curvature tensors.

    \begin{theorem}\label{thm:bundle_CY_intro}
        Let $E$ be a holomorphic vector bundle over a compact Hermitian manifold $(M, \om_g)$. Suppose that there exists a smooth Hermitian metric $h^E$ on $E$ such that
        \[
        \Lambda_{\om_g}\left(\sqrt{-1}\, R^{h^E}\right) > 0.
        \]
        Then, for any positive-definite Hermitian tensor $P \in \Gamma\left(M, E^* \otimes \bar{E}^*\right)$, there exists a unique smooth Hermitian metric $h$ on $E$ satisfying
        \beq
        \Lambda_{\om_g}\left(\sqrt{-1}\, R^{h}\right) = P.\label{PHYM0}
        \eeq
    \end{theorem}

\noindent  The prescribed Hermitian-Yang-Mills tensor equation \eqref{PHYM0} interpolates between the Calabi-Yau equation \eqref{CY} and the Hermitian-Einstein equation \eqref{HE}. By inheriting the geometric and analytic features of both, it serves as a pivotal bridge connecting the geometry of vector bundles to the broader realm of K\"ahler and Hermitian geometry.
To dynamically construct these canonical metrics, we introduce the  prescribed Hermitian--Yang--Mills flow, which serves as the parabolic analogue of the elliptic equation \eqref{PHYM0}:
    \begin{equation}
    \label{PHYMF}
    \begin{cases}
    \displaystyle \frac{\partial h}{\partial t} = -\Lambda_{\om_g}\left(\sqrt{-1}\, R^h\right) + P, \\[10pt]
    h(0) = h_0.
    \end{cases}\end{equation}
This flow is analogous to the Hermitian-Yang-Mills flow (\cite{Don85}), the Ricci flow (\cite{Ham82}), the K\"ahler-Ricci flow (\cite{Cao85}),   and  the Chern-Ricci flow (\cite{TW15}), providing a heat-flow method to deform an initial metric toward the  solution guaranteed by Theorem \ref{thm:bundle_CY_intro}. More generally, we consider a broader class of flows of the form
\begin{equation}
\label{GHYMF}
\frac{\partial h}{\partial t} = -\Lambda_{\omega_g}\left(\sqrt{-1}\, R^h\right) + \mathcal{F}_0(h),
\end{equation}
where $\mathcal{F}_0$ is a general map taking values in $\mathrm{Herm}(E)$; in particular, it includes the curvature associated with a Higgs bundle. Specific instances include:
\begin{enumerate}[label=(\roman*)]
    \item $\mathcal{F}_0(h) = P$ for some fixed $P \in \mathrm{Herm}(E)$;
    \item $\mathcal{F}_0(h) = \lambda h$ for $\lambda \in C^\infty(M, \R)$;
    \item $\mathcal{F}_0(h) = \lambda h + P + (\mathrm{tr}_{h_0} h) h_0$ for some fixed Hermitian metric $h_0$ on $E$.
\end{enumerate}

\noindent
It is well-known that  a priori uniform $C^0$-estimates are crucial for the  geometric flow theory. The first main result of this paper is a parabolic comparison principle, which yields uniform $C^0$-estimates for the flows under consideration. For simplicity, we present a basic version here that applies to the prescribed Hermitian--Yang--Mills flow \eqref{PHYMF}. A general version is established in Theorem~\ref{comparison}, providing uniform $C^0$-estimates for a broad class of flows of the form \eqref{GHYMF}.

    \btheorem\label{main2} Let $(M,\om_g)$ be a compact Hermitian manifold, and  $E$ be a holomorphic vector bundle on $M$. Suppose that $\underline h(t)$,   $h(t)$ and $\overline h(t)$ are smooth time-dependent Hermitian metrics on $E$ defined for $t\in [0,T)$ with $T\leq \infty$.  \bd \item[(1)]  If  $\overline h(0)\geq h(0)$ and
    \beq \frac{\p \overline h}{\p t}-\frac{\p h}{\p t}
    \geq \Lambda_{\omega_g}\left(\sqrt{-1} R^h\right)-\Lambda_{\omega_g}\left(\sqrt{-1} R^{\overline h}\right),\eeq
    then for any $t\in [0, T)$, one has \beq \overline h(t)\geq h(t).\eeq
    \item[(1')] If  $ \underline h(0)\leq h(0)$ and
    \beq \frac{\p \underline h}{\p t}-\frac{\p h}{\p t}
    \leq \Lambda_{\omega_g}\left(\sqrt{-1} R^h\right)-\Lambda_{\omega_g}\left(\sqrt{-1} R^{\underline h}\right),\eeq
    then for any $t\in [0, T)$, one has \beq  \underline h(t)\leq h(t).\eeq
    \ed

    \etheorem

\noindent As an application of Theorem \ref{main2}, we establish the following long-time convergence and existence results, which generalize Theorem \ref{thm:bundle_CY_intro}:

\btheorem\label{main3}Let $E$ be a holomorphic vector bundle over a compact Hermitian manifold $(M, \omega_g)$. Suppose that there exists a smooth Hermitian metric $h$ on $E$ such that $$
\Lambda_{\omega_g}\left(\sqrt{-1}\, R^{h}\right) $$ is quasi-positive as a Hermitian tensor in $\Gamma(M, E^* \otimes \bar{E}^*)$.
Then, for any quasi-positive tensor $P \in \Gamma(M, E^* \otimes \bar{E}^*)$ and any initial Hermitian metric $h_0$ on $E$, the prescribed Hermitian--Yang--Mills flow \eqref{PHYMF} exists for all $t\in [0,\infty)$ and converges smoothly at an exponential rate  to the \emph{unique} smooth Hermitian metric $h_\infty$ on $E$ satisfying
\begin{equation}
\Lambda_{\omega_g}\left(\sqrt{-1}\, R^{h_\infty}\right) = P.\label{PHYMinfty}
\end{equation}
\etheorem

\noindent We outline the proof of Theorem~\ref{main3}. First, we construct two smooth Hermitian barrier metrics $\underline{h}$ and $\overline{h}$ on $E$ satisfying
\begin{equation}
\underline{h} \leq h_0 \leq \overline{h}, \qquad
\Lambda_{\om_g}\left(\sqrt{-1}\, R^{\underline{h}}\right) \leq P \leq \Lambda_{\om_g}\left(\sqrt{-1}\, R^{\overline{h}}\right).
\end{equation}
By applying Theorem~\ref{main2}, we obtain the uniform $C^0$-estimate along the flow \eqref{PHYMF}:
\begin{equation}
\underline{h} \leq h(t) \leq \overline{h}.
\end{equation}
By using the evolution equation, we prove  the following exponential decay estimate:
\begin{equation}
\sup_M \left|\Lambda_{\omega_g}\left(\sqrt{-1}\, R^{h(t)}\right) - P\right|_{h(t)} \leq C_{2} e^{-C_1 t}.
\end{equation}
Combining these uniform bounds with the exponential decay, we conclude the long-time existence and convergence of the flow to a limiting metric $h_\infty$ satisfying \eqref{PHYMinfty}. A central innovation is that our proof entirely \emph{avoids} Donaldson's functional, relying solely on direct geometric estimates. Consequently, our results are valid on general Hermitian manifolds, unlike the classical theory which necessitates the K\"ahler condition.\\

 The following theorem extends the strategy established in the proof of Theorem \ref{main3}:
 \btheorem\label{main4}
  Let  $E$ be a holomorphic vector bundle over a compact Hermitian manifold $(M,\omega_g)$, and  $P\in\Gamma(M,E^*\ts\bar E^*)$ be a smooth Hermitian tensor. If there exist two smooth Hermitian metrics $\underline h$ and $\overline h$ on $E$ satisfying
 \beq
 \underline h\leq\overline h,
 \qquad
 \Lambda_{\omega_g}\left(\sq R^{\underline h}\right)
 \leq P\leq
 \Lambda_{\omega_g}\left(\sq R^{\overline h}\right), 
 \label{monotone barrier assumptions}
 \eeq
 then the prescribed Hermitian-Yang-Mills flows 
 \begin{equation}
 	\label{PHYMF2}
 	\begin{cases}
 		\displaystyle \frac{\partial h^-}{\partial t} = -\Lambda_{\om_g}\left(\sqrt{-1}\, R^{h^-}\right) + P, \\[10pt]
 		h^-(0) =\underline h, 
 \end{cases} \qquad \begin{cases}
 \displaystyle \frac{\partial h^+}{\partial t} = -\Lambda_{\om_g}\left(\sqrt{-1}\, R^{h^+}\right) + P, \\[10pt]
 h^+(0) =\overline h, 
\end{cases}  \end{equation}
 exist for all time $t\in[0,\infty)$. Moreover, \bd \item  $h^-(t)$ is non-decreasing, and $h^+(t)$ is non-increasing. 
 \item There exist smooth Hermitian metrics $h^-_\infty$ and $h^+_\infty$ on $E$ such that
 \beq
 h^-(t)\> h^-_\infty,
 \qquad
 h^+(t)\> h^+_\infty
 \quad\text{in }C^\infty,
 \eeq
 and $
 \underline h\leq h^-_\infty\leq h^+_\infty\leq\overline h$, 
 \beq
 \Lambda_{\omega_g}\left(\sq R^{h^-_\infty}\right)
 =P=
 \Lambda_{\omega_g}\left(\sq R^{h^+_\infty}\right).
 \eeq
\item  If $h$ is any smooth Hermitian metric  on $E$ satisfying $
 \underline h\leq h\leq\overline h$ and
 $$
 \Lambda_{\omega_g}\left(\sq R^h\right)=P,
 $$
 then
 \beq
 h^-_\infty\leq h\leq h^+_\infty.
 \eeq
 \ed 
 \etheorem

\noindent 
Note that if $P$ is not quasi-positive,  the solutions $h_\infty^-$ and $h_\infty^+$ in Theorem~\ref{main4} may differ, as can be seen clearly from the unitary-flat bundle case. This contrasts sharply​ with Theorem~\ref{main3}.
 The following result is a direct consequence of Theorem~\ref{main4}, and can be viewed as a \emph{compactness theorem}.
 
 \bcorollary\label{PHYM2}
 Let  $E$ be a holomorphic vector bundle over a compact Hermitian manifold $(M,\omega_g)$, and  $P\in\Gamma(M,E^*\ts\bar E^*)$ be a smooth Hermitian tensor. If there exist two smooth Hermitian metrics $\underline h$ and $\overline h$ on $E$ satisfying \beq
 \underline h\leq\overline h,
 \qquad
 \Lambda_{\omega_g}\left(\sq R^{\underline h}\right)
 \leq P\leq
 \Lambda_{\omega_g}\left(\sq R^{\overline h}\right),
 \eeq then there exists a smooth Hermitian metric $h$ such that
 \beq
 \underline h\leq h\leq\overline h,
 \qquad
 \Lambda_{\omega_g}\left(\sq R^h\right)=P.
 \label{PHYM3}
 \eeq
 \ecorollary

 \vskip 1\baselineskip

 In Theorem~\ref{main4}, the prescribed Hermitian-Yang-Mills flow converges globally for the specific initial data $\underline h$ and $\overline h$.  On a compact K\"ahler or balanced manifold, we obtain the following improvement. 
 
 \btheorem\label{main1}
Let  $E$ be a holomorphic vector bundle over a compact K\"ahler (or balanced) manifold
$(M,\omega_g)$.  For given Hermitian tensor $P\in\Gamma(M,E^*\ts \bar E^*)$,  if there exist two  Hermitian metrics \(\underline h\) and \(\overline h\) on $E$ such that
\beq
\underline h\leq  \overline h, \qtq{and}
\Lambda_{\omega_g}\left(\sq R^{\underline h}\right)\leq P\leq \Lambda_{\omega_g}\left(\sq R^{\overline
    h}\right),
\eeq
then for any  initial metric $h_0$ satisfying  $  \underline h\leq h_0\leq \overline h$,  the
prescribed Hermitian-Yang-Mills flow \eqref{PHYMF}
exists for all $t\in [0,\infty)$ and $ \underline h\leq h(t)\leq \overline h$. Moreover,  for any sequence $\{t_i\}\>\infty$, there exists a subsequence $\{t_{i_m}\}$ such that $\{h(t_{i_m})\}$ converges smoothly to a Hermitian metric $h_\infty$ on $E$ satisfying
\beq
\Lambda_{\omega_g}\left(\sq R^{h_\infty}\right)=P. \label{PHYM}
\eeq
\etheorem
\noindent
As we pointed out in Example \ref{example}, the prescribed Hermitian tensor $P$ in Theorem \ref{main1} may have  a negative
eigenvalue at some point, and \emph{it need not be positive}.  The proof of Theorem \ref{main1} proceeds as follows:  The uniform $C^0$-estimate along the flow is established by Theorem \ref{main2}; higher-order estimates and long-time existence are proved in Theorem \ref{parabolicbootstrap} and  Theorem \ref{thm long time existence}; finally, the long-time convergence is established via an analog of the Donaldson functional stated in Theorem \ref{thm modified Donaldson first variation}. Further generalizations of Theorem \ref{main3} and Theorem \ref{main1}, along with detailed discussions and examples, are presented in Section \ref{more}.

\bremark
\bd
\item By Serre duality, the corresponding statements of Theorem~\ref{main1} remain valid for negative bundles; see \cite{WYY26+}.  Building on the analogy with results established in \cite{FWYY26+}, we demonstrate that these principal results also hold for Higgs bundles over compact Hermitian manifolds.

\item  In Theorem \ref{comparison}, we established the uniform $C^0$-estimate for a  general class of flows \beq \frac{\partial h}{\partial t} = -\Lambda_{\omega_g}\left(\sqrt{-1}\, R^h\right) + \mathcal{F}_0(h). \eeq
The corresponding convergence theory will be addressed in future work.  For instance, similar results in this paper hold for the twisted Hermitian-Yang-Mills flow:
\beq \frac{\partial h}{\partial t} = -\Lambda_{\omega_g}\left(\sqrt{-1}\, R^h\right) +\lambda h+P, \eeq
upon straightforward modification to the arguments, as indicated in \cite{WYY26b+}.
\item Analogous results are expected to hold for the flow defined by RC-positivity and will be pursued in separate projects.
\ed
\eremark
    \noindent\textbf{Acknowledgements}. The second-named author would like to thank   Bing-Long Chen,  Jixiang Fu,  Kefeng Liu, Valentino Tosatti and Xi-Ping Zhu for inspiring  discussions.

\vskip 2\baselineskip
\section{Evolution equations and higher-order estimates}
Throughout this paper, all manifolds are assumed to be connected.
Let $(M,\omega_g) $ be a compact Hermitian manifold of dimension
$n$, and $E$ be a holomorphic vector bundle over $M$ of rank $r$.
Given a smooth Hermitian metric $h$ on $E$, the Chern connection of
$(E,h)$ is denoted by $\nabla^{h}$. We shall use the natural
decomposition $ \nabla^{h}=\p^{h}+\bp$, where $\p^{h}$ is the
$(1,0)$ part and $\bp$ is the $(0,1)$ part. Let $R^{h}$ be the Chern
curvature tensor of $(E,h)$. In local holomorphic coordinates
$\{z^i\}$ of $M$ and local holomorphic frame $\{e_\alpha\}$ of $E$,
\beq R^h=R_{i\bar j\alpha\bar\beta}dz^i\wedge d\bar z^j\ts
e^\alpha\ts\bar{e}^\beta \in\Gamma\left(M,\wedge^{1,1}T^*M\ts
E^*\ts\bar{E}^*\right), \eeq where \beq R_{i\bar j\alpha\bar\beta}
=  -\frac{\p^2 h_{\alpha\bar\beta}}{\p z^i\p \bar z^j}
+h^{\gamma\bar\delta} \frac{\p h_{\alpha\bar\delta}}{\p z^i}\frac{\p
h_{\gamma\bar\beta}}{\p \bar z^j}. \label{curvature tensor local
formula} \eeq The \textbf{Hermitian-Yang-Mills tensor} $ S^{h}\in
\Gamma(M,E^*\ts \bar E^*)$ of $(E,h)$  is defined as \beq S^{h}
:=\Lambda_{\omega_g}\left(\smo R^{h}\right)=\left(g^{i\bar
j}R_{i\bar j\alpha\bar\beta}\right) e^\alpha\ts \bar e^\beta.
\eeq To make the computations transparent, we adopt the following
notation throughout the paper.  By using the Hermitian $h$, there
are natural lifts of $R^h$ and $S^h$: \beq\Theta^h=R^h\cdot
h^{-1}=R_{i\bar j\alpha}^{\beta}dz^i\wedge d\bar z^j\ts
e^\alpha\ts{e}_\beta\in \Gamma(M,\Lambda^{1,1}T^*M\ts E^*\ts E)\eeq
and \beq   K^h=S^h\cdot h^{-1}=\left(g^{i\bar j}R_{i\bar
j\alpha}^{\beta}\right) e^\alpha\ts  e_\beta\in \Gamma(M,E^*\ts E).
\eeq
We say that $(M,\omega_g)$ is a \emph{Gauduchon manifold} if
$
\p\bp\om_g^{n-1}=0.
$
It is well known that, for any Hermitian metric $\omega$ on $M$,
there exists a smooth function $f$ such that
$
\omega_g = e^{f}\omega
$
satisfies $\p\bp\om_g^{n-1}=0$ (\cite{Gau84}). A Hermitian metric $\omega_g$ is called \emph{balanced} if $d\omega_g^{n-1}=0$, or equivalently $d^*\omega_g=0$.  We use the following notations:

\begin{itemize}
    \item $\mathrm{Herm}(E)$: the space of Hermitian tensors in
    $\Gamma(M, E^*\otimes \bar{E}^*)$;

    \item $\Herm^+(E)$: the subspace consisting of
    positive-definite Hermitian tensors;

    \item $\Herm^{\ge 0}(E)$: the subspace consisting of
    non-negative Hermitian tensors.
\end{itemize}
Recall that for $ P=P_\alpha^\beta e^\alpha\ts e_\beta \in \Omega^p(M,E^*\otimes E) $ and $ Q=Q_\alpha^\beta e^\alpha\ts e_\beta \in \Omega^q(M,E^*\otimes E) $, the product \beq  P\cdot Q \in \Omega^{p+q}(M,E^*\otimes E) \eeq  is defined as
\beq P \cdot Q = (P_\alpha^\beta \wedge Q_\beta^\gamma) \; e^\alpha \otimes e_\gamma. \eeq

\subsection{Some evolution equations} The prescribed Hermitian-Yang-Mills flow \eqref{PHYMF} can be written locally as
\beq
\frac{\p h_{\alpha\bar\beta}}{\p t} = -S^{h}_{\alpha\bar\beta}+P_{\alpha\bar\beta}.
\label{PHYM local lower index flow}
\eeq
The convergence of the term $$S^h-P$$ is the central focus of this paper. For notational simplicity, we define
\beq
\tilde S =S^h-P.
\eeq
We also use the notation
$$\Delta_{\C} f=\Lambda_{\omega_g}\left(\sq \p\bp f\right) $$  for $f\in C^\infty(M,\R)$.
\blemma\label{lem tilde S evolution}The evolution equation of $|\tilde S|^2_h$ along the flow \eqref{PHYMF} is
\beq
\left(\frac{\p }{\p t}-\Delta_\C \right)|\tilde S|_{h}^2
=-2\tilde S_{\alpha\bar\beta}  P_{\gamma\bar\delta}\tilde S_{\lambda\bar \mu} h^{\alpha\bar \mu}  h^{\lambda \bar \delta} h^{\gamma\bar\beta}
- \left|\nabla \tilde S\right|_{g,h}^2,
\label{tilde S evolution}
\eeq
where $\nabla$ denotes the connection on $E^*\ts \bar E^*$ induced by $(E,h)$.
\elemma

\bproof
The flow equation $\frac{\p h_{\alpha\bar\beta}}{\p t} =-\tilde S_{\alpha\bar\beta}$ gives
\beq
\frac{\p h^{\alpha\bar\beta}}{\p t}
= -h^{\alpha\bar \mu}\frac{\p h_{\lambda\bar \mu}}{\p t} h^{\lambda\bar\beta}
= h^{\alpha\bar \mu}\tilde S_{\lambda\bar \mu} h^{\lambda\bar\beta}.
\label{metric inverse evolution}
\eeq
Using the Christoffel symbols $\Gamma_{i\alpha}^{\gamma}=h^{\gamma\bar\delta}\frac{\p h_{\alpha\bar\delta}}{\p z^i}$, we have
\begin{eqnarray}
\frac{\p}{\p t}\Gamma_{i\alpha}^{\gamma}
&=&
\frac{\p h^{\gamma\bar\delta}}{\p t}
\frac{\p h_{\alpha\bar\delta}}{\p z^i}
+
h^{\gamma\bar\delta}
\frac{\p}{\p z^i}\left(\frac{\p h_{\alpha\bar\delta}}{\p t}\right)\nonumber\\
&=&
h^{\gamma\bar\eta}\tilde S_{\lambda\bar\eta}h^{\lambda\bar\delta}
\frac{\p h_{\alpha\bar\delta}}{\p z^i}
-
h^{\gamma\bar\delta}
\frac{\p \tilde S_{\alpha\bar\delta}}{\p z^i}=
-h^{\gamma\bar\delta}\nabla_i\tilde S_{\alpha\bar\delta}.
\label{Gamma evolution}
\end{eqnarray}
The Chern curvature formula $$R^{h}_{i\bar j\alpha\bar\beta}=-\frac{\p \Gamma_{i\alpha}^{\gamma}}{\p \bar z^j} h_{\gamma\bar\beta}$$  yields
\begin{eqnarray}
\frac{\p}{\p t}R^{h}_{i\bar j\alpha\bar\beta}&=& - \frac{\p}{\p \bar z^j} \left(\frac{\p \Gamma_{i\alpha}^{\gamma}}{\p t}\right) h_{\gamma\bar\beta}
-  \frac{\p \Gamma_{i\alpha}^{\gamma}}{\p \bar z^j} \frac{\p h_{\gamma\bar\beta}}{\p t} \nonumber\\
&=&\frac{\p}{\p \bar z^j}
\left(h^{\gamma\bar\delta}\nabla_i\tilde S_{\alpha\bar\delta}\right)
h_{\gamma\bar\beta}
+
\frac{\p \Gamma_{i\alpha}^{\gamma}}{\p \bar z^j}
\tilde S_{\gamma\bar\beta}\nonumber\\
&=&\nabla_{\bar j}\nabla_i\tilde S_{\alpha\bar\beta}
-
R^{h}_{i\bar j\alpha\bar\delta}
h^{\gamma\bar\delta}\tilde S_{\gamma\bar\beta}.
\end{eqnarray}
In particular, we deduce that
\begin{eqnarray}
\frac{\p}{\p t}\tilde S_{\alpha\bar\beta}
=
\frac{\p}{\p t}S^{h}_{\alpha\bar\beta}
= g^{i\bar j}\frac{\p}{\p t}R^{h}_{i\bar j\alpha\bar\beta} &=&
g^{i\bar j}\nabla_{\bar j}\nabla_i\tilde S_{\alpha\bar\beta}
-
S^{h}_{\alpha\bar\delta}
h^{\gamma\bar\delta}\tilde S_{\gamma\bar\beta}\nonumber\\
&=&
g^{i\bar j}\nabla_i\nabla_{\bar j}\tilde S_{\alpha\bar\beta}
-
\tilde S_{\alpha\bar\gamma}
h^{\delta\bar\gamma}S^{h}_{\delta\bar\beta}. \label{evolutionS}
\end{eqnarray}
Here  the last equality follows from the Ricci identity
\beq
\nabla_{\bar j}\nabla_i\tilde S_{\alpha\bar\beta}
=
\nabla_i\nabla_{\bar j}\tilde S_{\alpha\bar\beta}
+
R^{h}_{i\bar j\alpha\bar\delta}
h^{\gamma\bar\delta}\tilde S_{\gamma\bar\beta}
-
R^{h}_{i\bar j\delta\bar\beta}
h^{\delta\bar\gamma}\tilde S_{\alpha\bar\gamma}.
\eeq
Since $|\tilde S|_{h}^2=\tilde S_{\alpha\bar\beta}\tilde S_{\gamma\bar\delta} h^{\alpha\bar\delta}h^{\gamma\bar\beta}$ and $\tilde S$ is Hermitian,
\begin{eqnarray}
\frac{\p}{\p t}|\tilde S|_{h}^2
&=&  \frac{\p \tilde S_{\alpha\bar\beta}}{\p t} \tilde S_{\gamma\bar\delta} h^{\alpha\bar\delta}h^{\gamma\bar\beta}
+ \tilde S_{\alpha\bar\beta}  \frac{\p \tilde S_{\gamma\bar\delta}}{\p t}  h^{\alpha\bar\delta}h^{\gamma\bar\beta} \nonumber\\
&&+  \tilde S_{\alpha\bar\beta}\tilde S_{\gamma\bar\delta} h^{\alpha\bar \mu} \tilde S_{\lambda\bar \mu} h^{\lambda\bar \delta} h^{\gamma\bar\beta}
+ \tilde S_{\alpha\bar\beta}\tilde S_{\gamma\bar\delta} h^{\alpha\bar\delta} h^{\gamma\bar \mu}\tilde S_{\lambda\bar \mu} h^{\lambda\bar\beta}\nonumber\\
&=& \LL \p_t \tilde S,\tilde S\RL_{h} +\LL \tilde S,\p_t \tilde S\RL_{h} +
2 \tilde S_{\alpha\bar\beta}\tilde S_{\gamma\bar\delta}\tilde S_{\lambda\bar \mu} h^{\alpha\bar \mu}  h^{\lambda\bar \delta} h^{\gamma\bar\beta}.
\label{time derivative of tilde S norm}
\end{eqnarray}
 The Ricci identity yields
\begin{eqnarray}
g^{i\bar j}\nabla_{\bar j}\nabla_{i}\tilde S_{\alpha\bar\beta}
- g^{i\bar j}\nabla_i\nabla_{\bar j}\tilde S_{\alpha\bar\beta}
&=&  g^{i\bar j}R^{h}_{i\bar j\alpha\bar\delta}h^{\gamma\bar\delta} \tilde S_{\gamma\bar\beta}
-g^{i\bar j}R^{h}_{i\bar j\gamma\bar\beta}h^{\gamma\bar\delta} \tilde S_{\alpha\bar\delta}\nonumber\\
&=& S^{h}_{\alpha\bar\delta}h^{\gamma\bar\delta} \tilde S_{\gamma\bar\beta}
- S^{h}_{\gamma\bar\beta}h^{\gamma\bar\delta} \tilde S_{\alpha\bar\delta}\nonumber\\
&=& - S^{h}_{\alpha\bar\delta}h^{\gamma\bar\delta} P_{\gamma\bar\beta}
+S^{h}_{\gamma\bar\beta}h^{\gamma\bar\delta} P_{\alpha\bar\delta}.
\label{Ricci identity for tilde S}
\end{eqnarray}
A straightforward computation shows
\beq
\LL g^{i\bar j}\nabla_{\bar j}\nabla_{i}\tilde S
- g^{i\bar j}\nabla_i\nabla_{\bar j}\tilde S,\tilde S\RL_h=0.
\eeq
Therefore,
\begin{eqnarray}
\Delta_{\C}|\tilde S|_{h}^2
&=&   \LL g^{i\bar j}\nabla_i\nabla_{\bar j}\tilde S,\tilde S\RL_{h}
+\LL \tilde S,g^{j\bar i}\nabla_{\bar i}\nabla_j\tilde S\RL_{h} + \left|\nabla \tilde S\right|_{g,h}^2\nonumber\\
&=&   2\mathrm{Re}\left(\LL g^{i\bar j}\nabla_i\nabla_{\bar j}\tilde S,\tilde S\RL_{h}\right)
+ \left|\nabla \tilde S\right|_{g,h}^2.
\label{Laplacian tilde S norm}
\end{eqnarray}
Combining \eqref{time derivative of tilde S norm} and \eqref{Laplacian tilde S norm}, we get
\beq
\left(\frac{\p }{\p t}-\Delta_\C \right)|\tilde S|_{h}^2
= 2\mathrm{Re}\left(\LL \p_t \tilde S -  g^{i\bar j}\nabla_i\nabla_{\bar j}\tilde S,\tilde S\RL_{h}\right)
+ 2 \tilde S_{\alpha\bar\beta}\tilde S_{\gamma\bar\delta}\tilde S_{\lambda\bar \mu} h^{\alpha\bar \mu}  h^{\lambda\bar \delta} h^{\gamma\bar\beta}
- \left|\nabla \tilde S\right|_{g,h}^2 .
\label{pre tilde S evolution}
\eeq
By the evolution equation \eqref{evolutionS}, one has
\beq
2\mathrm{Re}\left(\LL \p_t \tilde S -  g^{i\bar j}\nabla_i\nabla_{\bar j}\tilde S,\tilde S\RL_h \right)
= -2\tilde S_{\alpha\bar\beta}  S^{h}_{\gamma\bar\delta}\tilde S_{\lambda\bar \mu} h^{\alpha\bar \mu}  h^{\lambda\bar \delta} h^{\gamma\bar\beta}.
\label{tilde S time minus laplacian term}
\eeq
Substituting \eqref{tilde S time minus laplacian term} into \eqref{pre tilde S evolution} yields
\begin{eqnarray}
\left(\frac{\p }{\p t}-\Delta_\C \right)|\tilde S|_h^2
&=& -2\tilde S_{\alpha\bar\beta}  S^{h}_{\gamma\bar\delta}\tilde S_{\lambda\bar \mu} h^{\alpha\bar \mu}  h^{\lambda\bar \delta} h^{\gamma\bar\beta}
+ 2 \tilde S_{\alpha\bar\beta}\tilde S_{\gamma\bar\delta}\tilde S_{\lambda\bar \mu} h^{\alpha\bar \mu}  h^{\lambda\bar \delta} h^{\gamma\bar\beta}
- \left|\nabla \tilde S\right|_{g,h}^2  \nonumber\\
&=& -2\tilde S_{\alpha\bar\beta}  P_{\gamma\bar\delta}\tilde S_{\lambda\bar \mu} h^{\alpha\bar \mu}  h^{\lambda\bar \delta} h^{\gamma\bar\beta}
- \left|\nabla \tilde S\right|_{g,h}^2.
\end{eqnarray}
This is \eqref{tilde S evolution}.
\eproof

\bcorollary\label{cor tilde S inequality} The following estimate holds:
\beq
\left(\frac{\p }{\p t}-\Delta_\C \right)|\tilde S|_{h}^2 \leq -2 \lambda^P_{\min}\cdot |\tilde S|_{h}^2,
\eeq
where
\beq \lambda^P_{\min}(x,t):=\inf_{v \in E_x, v\neq 0}  \frac{ P(v, \bar v) }{h(t)(v,\bar v)}.\eeq

\ecorollary

\bproof Fix a point $q\in M$ and a time $t$.  Choose an $h(t)$-unitary basis of
$E$ such that
$$h_{\alpha\bar\beta}=\delta_{\alpha\beta}, \qtq{and} P_{\alpha\bar\beta}=\lambda_\alpha\delta_{\alpha\beta}, $$ where
$\lambda_1\leq\cdots\leq\lambda_r$. One has
\[
\tilde S_{\alpha\bar\beta}P_{\gamma\bar\delta}\tilde S_{\lambda\bar\mu}
h^{\alpha\bar\mu}h^{\lambda\bar\delta}h^{\gamma\bar\beta}=\sum_{\alpha,\beta}\lambda_\beta|\tilde S_{\alpha\bar\beta}|^2
\geq \lambda_1\sum_{\alpha,\beta}|\tilde S_{\alpha\bar\beta}|^2
=\lambda^P_{\min}|\tilde S|_h^2.
\]
The desired estimate now follows from Lemma~\ref{lem tilde S evolution} by discarding $\left|\nabla \tilde S\right|_{g,h}^2$.
\eproof

\noindent The following transformation formula for the
Hermitian-Yang-Mills curvature tensors is essentially well-known
(e.g. \cite[Proposition~2.8]{FWYY26+}) \blemma Let $(M,\omega_g)$ be
a compact  Hermitian manifold and $E\>M$ be a holomorphic vector
bundle. If $h$ and $h_0$ are smooth Hermitian metrics on $E$, then
\beq \Lambda_{\omega_g}\left(\sq  \Theta^{h}\right)-
\Lambda_{\omega_g}\left(\sq  \Theta^{h_0}\right)=
\Lambda_{\omega_g}\left(\smo\bp\left( (\p^{h_0} H) \cdot H^{-1}
\right)\right)\label{conformalchange1} \eeq as tensors in $\Gamma(M,E^*\ts
E)$ and $H=h\cdot h_0^{-1}$.

\elemma

\subsection{The long-time existence under uniform $C^0$-estimates} In this subsection, we establish a long-time existence theorem for the flow \eqref{PHYMF} on general Hermitian manifolds, under certain assumptions.

\btheorem\label{thm finite time extension}
Let $(M,\omega_g)$ be a compact Hermitian manifold, and $E$ be a holomorphic vector bundle over $M$.
Suppose that  $h_t$, $0\leq t<T<\infty$, is a smooth solution to the prescribed Hermitian-Yang-Mills flow \eqref{PHYMF}. If there exist constants $C_1$ and $C_2$ such that
\beq
C_1^{-1}h_0\leq h_t\leq C_1h_0,
\qquad
\sup_M\left|S^{h_t}-P\right|_{h_t}\leq C_2
\label{finite extension assumptions}
\eeq
for all $t\in[0,T)$, then $\{h_t\}$ converges in $C^\infty$ to a smooth Hermitian metric $h_T$ as $t\to T$. In particular, the flow extends smoothly beyond $T$.
\etheorem
\noindent
The proof of Theorem \ref{thm finite time extension} relies on a careful adaptation of the classical theory of elliptic and parabolic equations. While the existing literature primarily addresses the K\"ahler case with $P = 0$ (e.g., \cite{Don85}, \cite{Siu87}), we provide an independent proof applicable to the non-K\"ahler setting, thereby ensuring the exposition remains self-contained.\\

    Let $\mathrm{Herm}^+(E)$ be the space of Hermitian metrics on $E$.  The right-hand side of \eqref{PHYMF} defines a smooth second-order differential operator
\[
L:\mathrm{Herm}^+(E)\longrightarrow \mathrm{Herm}(E).
\]
For $h\in\mathrm{Herm}^+(E)$ and $\eta\in\mathrm{Herm}(E)$, the linearization of $L$ is given by
\beq
DL_h(\eta)_{\alpha\bar\beta}=g^{i\bar j}\frac{\p^2\eta_{\alpha\bar\beta}}{\p z^i\p\bar z^j}+\text{lower order terms},
\eeq
which is strongly elliptic.  By appealing to the short-time existence theory for quasilinear parabolic systems on compact manifolds (e.g., \cite[Chapter~XIII, Section~7]{Lie96}), we obtain a unique smooth solution to \eqref{PHYMF} defined on a maximal time interval $[0, T)$ with $T>0$.
For higher-order estimates, we adopt the following conventions.   On the Hermitian vector bundle $(E,h_0)$ over $(M,\omega_g)$,  we define the $C^m$-norm of a smooth section $\phi\in\Gamma\left(M,E^*\ts\bar E^*\right)$ with respect to $\omega_g$ and $h_0$ by
\beq
\left\|\phi\right\|_{C^m(M,\omega_g, h_0)}
:=\sum_{j=0}^{m}\sup_M\left|\left(\nabla^{h_0}\right)^j \phi\right|_{g,h_0},
\label{global Ck norm definition}
\eeq
where $\nabla^{h_0}$ is the covariant derivative induced by Chern connections of $(E,h_0)$ and $(T^{1,0}M,\omega_g)$.  We fix a finite open cover $\mathcal{U} = \{U_\mu\}_{\mu=1}^N$ of $M$ such that each $U_\mu$ is relatively compact in a larger domain $\widetilde{U}_\mu$, equipped with holomorphic coordinates and a local holomorphic frame $\{e_1^{(\mu)}, \dots, e_r^{(\mu)}\}$ of $E$.
With respect to these local frames, the global $C^m$-norm \eqref{global Ck norm definition} is equivalent to the local coordinate norm
\begin{equation}
\sum_{\mu} \sum_{\alpha,\beta} \sum_{|I| \le m}
\sup_{U_\mu} \bigl| \partial^I \phi_{\alpha\bar\beta} \bigr|,
\label{local coordinate Ck norm}
\end{equation}
where $I$ denotes a multi-index. The constants implicit in this equivalence depend only on $\omega_g$, $h_0$, and the chosen finite cover. We also employ parabolic norms locally, after fixing the coordinate charts and frames described above. For a cylinder $$Q=B_r(x_\mu)\times I\Subset U_\mu\times\R,$$ we write $x=(x^1,\cdots,x^{2n})$ for the associated real coordinates. For a scalar function $f(x,t)$ defined on $Q$ and $1<p<\infty$, we define the parabolic $W^{2,p}$-norm  by
\beq
\left\|f\right\|_{W^{2,p}_t(Q)}
:=\sum_{|\nu|\leq2}\left\|\p_x^\nu f\right\|_{L^p(Q)}
+\left\|\p_t f\right\|_{L^p(Q)}.
\label{parabolic W21p norm definition}
\eeq
For $0<\alpha<1$, define
\beq
d((x,t),(y,s)):=|x-y|+|t-s|^{1/2},
\qquad
[f]_{\alpha;Q}:=\sup_{\substack{(x,t),(y,s)\in Q\\(x,t)\neq(y,s)}}
\frac{|f(x,t)-f(y,s)|}{d((x,t),(y,s))^\alpha}.
\label{parabolic Holder seminorm definition}
\eeq
For an integer $k\geq0$, the  parabolic H\"older norm (e.g., \cite[Chapter~IV, Section~1]{Lie96}) is
\beq
\left\|f\right\|_{C^{k+\alpha,(k+\alpha)/2}(Q)}
:=\sum_{|\nu|+2j\leq k}\left\|\p_x^\nu\p_t^j f\right\|_{C^0(Q)}
+\sum_{|\nu|+2j=k}\left[\p_x^\nu\p_t^j f\right]_{\alpha;Q}.
\label{parabolic Holder norm definition}
\eeq

\noindent A crucial ingredient in the proof of Theorem \ref{thm finite time extension} is the higher-order estimate for the Hermitian metrics $h_t$ established below.

\btheorem\label{parabolicbootstrap}
Let $(M,\omega_g)$ be a compact Hermitian manifold, and $E$ be a holomorphic vector bundle over $M$.
Suppose that  $h_t$, $0\leq t<T\leq \infty$, is a smooth solution to the flow \eqref{PHYMF}. If there exist constants $C_1$ and $C_2$ such that
\beq
C_1^{-1}h_0\leq h_t\leq C_1h_0,
\qquad
\sup_M\left|S^{h_t}-P\right|_{h_t}\leq C_2
\eeq
for all $t\in[0,T)$, then for every $\tau\in (0, T)$ and every integer $m\geq 0$, there exists a positive constant $C_{m,\tau}=C_{m,\tau}(\omega_g,h_0,P,C_1,C_2)$ such that
\beq
\left\|h_t\right\|_{C^{m}(M,\omega_g, h_0)}\leq C_{m,\tau},
\qquad t\in[\tau,T).
\label{uniform parabolic bootstrapping estimate}
\eeq
\etheorem

\bproof We first show that  if a smooth Hermitian metric $k$ satisfies
\beq
C_1^{-1}h_0\leq k\leq C_1h_0,
\qquad
\sup_M\left|S^{k}\right|_{k}\leq C_2,
\label{C1 estimate assumption}
\eeq
then $\left\|k\right\|_{C^1(M,\omega_g, h_0)}$ is uniformly bounded by a constant depending only on $\om_g$, $h_0$, $C_1$ and $C_2$. Otherwise, there are Hermitian metrics $k_i$ satisfying \eqref{C1 estimate assumption} and points $x_i\to x_\infty$ such that
\beq
m_i:=\sup_M\left|\nabla^{h_0}k_i\right|_{h_0}
=\left|\nabla^{h_0}k_i\right|_{h_0}(x_i)\to+\infty.
\label{C1 blow up sequence}
\eeq
We choose holomorphic coordinates $\{z^A\}$ and a holomorphic frame $\{e_\alpha\}$ centered at $x_\infty$ such that $g_{A\bar B}(x_\infty)=\delta_{AB}$ and $h_{0,\alpha\bar\beta}(x_\infty)=\delta_{\alpha\beta}$.  We select $\rho>0$ so that the coordinate chart contains $D_{3\rho}$. After passing to a subsequence, we have $z_i=z(x_i)\in D_\rho$.
On $D_{2\rho}$, the Euclidean component norms associated with the chosen holomorphic coordinates and local frame are uniformly equivalent to the norms induced by $\omega_g$ and $h_0$.
Moreover, $$dk_i-\nabla^{h_0}k_i$$ is a zeroth-order expression in $k_i$.
Since the sequence $k_i$ is uniformly bounded, there exist constants $L_0\geq1$ and $L_1>0$, depending only on $\omega_g,h_0,C_1$ and the fixed finite system of coordinate charts and local frames, such that
\beq
L_0^{-1}m_i-L_1\leq |dk_i|(x_i),\qquad
\sup_{D_{2\rho}}|dk_i|\leq L_0m_i+L_1.
\label{ordinary derivative comparison in C1 proof}
\eeq
Here, $dk_i$ denotes the Euclidean gradient of the component matrix $k_{i,\alpha\bar{\beta}}$ with respect to the spatial coordinates, and its norm is computed via the standard Euclidean norm over all coordinate and bundle indices. Define rescaled tensors $\tilde k_i\in C^\infty\left(D_{m_i\rho}, \C^{r\times r}\right)$ by \beq \tilde k_i(w):=k_i(z_i+w/m_i).\eeq Given any $R > 0$, $\tilde{k}_i$ is defined on $D_R$ provided that $i$ is sufficiently large, and
\beq
|d\tilde k_i|(0)\geq \frac{1}{2L_0},
\qquad
\sup_{D_R}|d\tilde k_i|\leq 2L_0.
\label{normalized derivative in C1 proof}
\eeq
In the chosen frames,
\beq
S^{k_i}_{\alpha\bar\beta}
=-g^{A\bar B}\frac{\p^2 k_{i,\alpha\bar\beta}}{\p z^A\p\bar z^B}
+g^{A\bar B}k_i^{\gamma\bar\delta}
\frac{\p k_{i,\alpha\bar\delta}}{\p z^A}
\frac{\p k_{i,\gamma\bar\beta}}{\p\bar z^B}.
\label{local elliptic equation in C1 proof}
\eeq
Hence the rescaled metrics solve
\beq
g_i^{A\bar B}\frac{\p^2\tilde k_{i,\alpha\bar\beta}}{\p w^A\p\bar w^B}
=g_i^{A\bar B}\tilde k_i^{\gamma\bar\delta}
\frac{\p\tilde k_{i,\alpha\bar\delta}}{\p w^A}
\frac{\p\tilde k_{i,\gamma\bar\beta}}{\p\bar w^B}
-\frac{1}{m_i^2}S_{i,\alpha\bar\beta},
\label{rescaled elliptic equation in C1 proof}
\eeq
where $g_i^{A\bar B}(w)=g^{A\bar B}(z_i+w/m_i)$ and $S_i(w)=S^{k_i}(z_i+w/m_i)$.
By \eqref{normalized derivative in C1 proof} and \eqref{C1 estimate assumption}, the right-hand side of \eqref{rescaled elliptic equation in C1 proof} is uniformly bounded on each fixed $D_R$. Fix $p>2n$ and choose $\alpha\in (0, 1-2n/p)$. Applying the interior elliptic $W^{2,p}$-estimates on $D_{R+1}$ and the compact embedding $W^{2,p}(D_R)\hookrightarrow C^{1,\alpha}(\bar D_R)$, and then using a diagonal argument over $R\in\mathbb N$, we obtain a single subsequence $\tilde k_{i_s}$ and $\tilde k_\infty\in C^{1,\alpha}_{\mathrm{loc}}(\C^n,\C^{r\times r})$ such that, for every $R>0$,
\beq
\lim_{s\to\infty}\|\tilde k_{i_s}-\tilde k_\infty\|_{C^{1,\alpha}(\bar{D_R}, g_{\C^n})}=0, \qtq{and} \tilde k_{i_s}\rightharpoonup\tilde k_\infty
\quad\text{weakly in }W^{2,p}({D_R}, g_{\C^n}).
\label{C1 local convergence in C1 proof}
\eeq
In particular, $|d\tilde k_\infty|(0)\geq 1/(2L_0)$. Since $h_0(z_i+w/m_i)\to I_r$ locally uniformly in the chosen frame, the metric equivalence in \eqref{C1 estimate assumption} passes to the limit and gives
\begin{equation}
C_1^{-1}I_r\leq \tilde k_\infty\leq C_1I_r
\qquad\text{on }\C^n.
\label{uniform bounds for flat limit in C1 proof}
\end{equation}
By passing to the weak limit in \eqref{rescaled elliptic equation in C1 proof}, we deduce that $\tilde k_\infty$  satisfies the following elliptic equation on $\C^n$ in the distribution sense
\beq
\Delta_0 \tilde k_\infty=\sum_A(\p_A\tilde k_\infty)\tilde k_\infty^{-1}(\p_{\bar A}\tilde k_\infty),
\label{limiting flat elliptic equation in C1 proof}
\eeq
where $\Delta_0=\sum_A\frac{\p^2}{\p w^A \p\bar w^A}$.
Since the right-hand side of \eqref{limiting flat elliptic equation in C1 proof} is locally $C^\alpha$, elliptic bootstrapping shows that $\tilde k_\infty$ is smooth. In the following, we show that $\tilde k_\infty$ is constant. This will contradict $|d\tilde k_\infty|(0)\geq 1/(2L_0)$ and complete the proof of the asserted $C^1$-estimate.\\

Set $K=\tilde k_\infty$. Since $K$ is Hermitian and positive definite,
\[
(\p_AK)^*=\p_{\bar A}(K^*)=\p_{\bar A}K,
\qquad
(K^{-1/2})^*=K^{-1/2}.
\]
Using these identities, we define the non-negative energy density
\begin{eqnarray}
e
&:=&\sum_A\tr\left(
K^{-1}(\p_AK)K^{-1}(\p_{\bar A}K)
\right)\nonumber\\
&=&\sum_A\tr\left[
\left(K^{-1/2}(\p_AK)K^{-1/2}\right)
\left(K^{-1/2}(\p_AK)K^{-1/2}\right)^*
\right]\geq0.
\end{eqnarray}
If $e\equiv0$, then every summand vanishes. Since $K^{-1/2}$ is invertible, this gives $\p_AK=0$ for every $A$; the adjoint identity also gives $\p_{\bar A}K=0$, and so $K$ is constant. Thus, it remains to prove that $e\equiv0$.
Using \eqref{limiting flat elliptic equation in C1 proof}, a straightforward calculation shows
\be
\Delta_0e
&=&\sum_{A,B}\tr\left\{
K^{-1}\left(\p_A\p_BK-(\p_AK)K^{-1}(\p_BK)\right)
K^{-1}\left(\p_A\p_BK-(\p_AK)K^{-1}(\p_BK)\right)^*
\right\}\\
&&+\sum_{A,B}\tr\left\{
K^{-1}\left(\p_A\p_{\bar B}K-(\p_AK)K^{-1}(\p_{\bar B}K)\right)
K^{-1}\left(\p_A\p_{\bar B}K-(\p_AK)K^{-1}(\p_{\bar B}K)\right)^*
\right\}.
\ee
Since for any $r\times r$ matrix $A$,
\[
    \tr\left(K^{-1}AK^{-1}A^*\right)
    = \tr\left[\left(K^{-1/2}AK^{-1/2}\right)\left(K^{-1/2}AK^{-1/2}\right)^*\right]\geq 0,
\]
we conclude that 
\beq
    \Delta_0 e\geq 0.
\eeq 
Thus, $e$ is subharmonic. On the other hand, by \eqref{uniform bounds for flat limit in C1 proof} and \eqref{limiting flat elliptic equation in C1 proof},
\beq
e
=\sum_A\tr\left(
K^{-1}(\p_AK)K^{-1}(\p_{\bar A}K)
\right)
\leq C_1\sum_A\tr\left(
(\p_AK)K^{-1}(\p_{\bar A}K)
\right)
=C_1\Delta_0\tr K.
\eeq
Fix some $w_0\in\C^n$, and let $dV_0$ be the Euclidean volume form. There exists a constant $c_n>0$, depending only on $n$, such that for every $R>0$ one can choose a standard cutoff function $\eta_R$ satisfying
\beq
0\leq\eta_R\leq1,\qquad
\eta_R=1\ \text{on }B_R(w_0),\qquad
\supp\eta_R\subset B_{2R}(w_0),
\qquad
|\Delta_0\eta_R|\leq c_nR^{-2}.
\eeq
Integration by parts and $\tr K\leq rC_1$ give
\begin{eqnarray}
C_1^{-1}\int_{B_R(w_0)}e\,dV_0
&\leq&
\int_{\C^n}\eta_R\Delta_0\tr K\,dV_0
=\int_{\C^n}\tr K\,\Delta_0\eta_R\,dV_0\nonumber\\
&\leq&
\left|\int_{\C^n}\tr K\,\Delta_0\eta_R\,dV_0\right|
\leq c_n'R^{2n-2},
\end{eqnarray}
where $c_n'=c_n'(n,r,C_1)$ is a constant.
Since $4\Delta_0$ is the ordinary Euclidean Laplacian on $\C^n$, the mean-value inequality for the non-negative subharmonic function $e$ yields
\beq
e(w_0)
\leq
\frac{1}{|B_R(w_0)|}
\int_{B_R(w_0)}e\,dV_0
\leq \frac{c_n''}{R^2},
\eeq
where $|B_R(w_0)|=\int_{B_R(w_0)}dV_0$ and $c_n''=c_n''(n,r,C_1)$ is a constant.
Letting $R\to\infty$ and using the arbitrariness of $w_0$, we obtain $e\equiv0$. This completes the proof of the $C^1$-estimate under \eqref{C1 estimate assumption}.

\vone For $\{h_t\}$ along the flow \eqref{PHYMF}, one has \beq
\sup_M\left|S^{h_t}\right|_{h_t} \leq
\sup_M\left|S^{h_t}-P\right|_{h_t} +\sup_M|P|_{h_t} \leq
C_2+C_1\sup_M|P|_{h_0}. \label{S bound from tilde S bound} \eeq
Hence, there is a constant $C_3=C_3(\om_g,h_0,P,C_1,C_2)$ such that
\beq \left\|h_t\right\|_{C^1(M,\omega_g, h_0)}\leq C_3, \qquad t\in[0,T).
\label{lower parabolic input estimates} \eeq

For higher-order estimates, we fix $\tau\in(0,T)$,
$\alpha\in(0,1)$, and $p>\frac{2n+2}{1-\alpha}$. Fix
$\rho=\frac{\sqrt{\tau}}{4}$ with $4\rho^2=\tau/4<\tau$. Diminishing $\rho$ if necessary,  we  choose
finitely many holomorphic coordinate charts $U_\mu$ with local
holomorphic frames such that $B_{4\rho}(x_\mu)\Subset U_\mu$ and
$M=\bigcup_{\mu=1}^N B_{\rho/4}(x_\mu)$.  For $t_0\in[\tau,T)$, we
set \beq Q^{\mu}_{r}(t_0):=B_r(x_\mu)\times(t_0-r^2,t_0].
\label{local parabolic cylinders} \eeq Then
$Q^{\mu}_{2\rho}(t_0)\Subset U_\mu\times(0,T)$. In the coordinate
chart $U_\mu$, the flow \eqref{PHYMF} takes the form \beq
(\p_t-\Delta_\C)h_{\alpha\bar\beta}=F_{\alpha\bar\beta}, \qquad
F_{\alpha\bar\beta}=-g^{i\bar j}h^{\gamma\bar\delta} \frac{\p
h_{\alpha\bar\delta}}{\p z^i} \frac{\p h_{\gamma\bar\beta}}{\p\bar
z^j}+P_{\alpha\bar\beta}. \label{fixed scalar operator bootstrap}
\eeq The operator $\p_t-\Delta_\C$ is fixed and uniformly parabolic.
The uniform $C^1$-estimate \eqref{lower parabolic input estimates} gives a uniform
constant $C_4$ such that, for every $\alpha,\beta,\mu$ and every
$t_0\in[\tau,T)$, \beq
\left\|F_{\alpha\bar\beta}\right\|_{L^\infty(Q^\mu_{2\rho}(t_0))}\leq
C_4. \label{F Linfty bound} \eeq All constants introduced
henceforth, starting with $C_4$, depend (at most) on $\tau$,
$\omega_g$, $h_0$, $p$, $P$, $C_1$, and $C_2$; additionally,
$C_{m,\tau}$ may depend on $m$. Crucially, none of these constants
depend on $t_0$ or $T$, and this independence will not be stated
explicitly again. The interior parabolic $W^{2,p}_t$ estimate (e.g.,
\cite[Theorem~7.22]{Lie96}), applied on $Q^\mu_\rho(t_0)\subset
Q^\mu_{2\rho}(t_0)$, yields \beq
\left\|h_{\alpha\bar\beta}\right\|_{W^{2,p}_t(Q^\mu_\rho(t_0))}\leq
C_5. \label{W21p first estimate} \eeq The parabolic Sobolev
embedding theorem (e.g., \cite[Chapter~VI, Section~3]{Lie96}) gives
a uniform $C^{1+\alpha,(1+\alpha)/2}$-bound. Moreover, the expression in
\eqref{fixed scalar operator bootstrap} shows that
\beq \left\|F_{\alpha\bar\beta}\right\|_{C^{\alpha,\alpha/2}(Q^\mu_\rho(t_0))}\leq C_6.  \eeq  The
interior parabolic Schauder estimate (e.g.,
\cite[Theorem~4.9]{Lie96}) gives \beq
\left\|h_{\alpha\bar\beta}\right\|_{C^{2+\alpha,(2+\alpha)/2}\left(Q^\mu_{3\rho/4}(t_0)\right)}\leq
C_7. \label{first Schauder bootstrap} \eeq Repeating the argument on
nested cylinders between $Q^\mu_\rho(t_0)$ and $Q^\mu_{\rho/4}(t_0)$
yields \beq
\left\|h_{\alpha\bar\beta}\right\|_{C^{m+\alpha,(m+\alpha)/2}\left(Q^\mu_{\rho/4}(t_0)\right)}\leq
C_{m,\tau} \label{local high order parabolic estimate} \eeq  for
every $m \geq 0$. Setting $t = t_0$ and using the finite covering of
$M$, together with the equivalence between \eqref{global Ck norm
definition} and \eqref{local coordinate Ck norm}, we obtain
\eqref{uniform parabolic bootstrapping estimate}.
\eproof

\bproof[Proof of Theorem \ref{thm finite time extension} ]
Taking $\tau=T/2$ in Theorem \ref{parabolicbootstrap}, for every integer $m\geq0$ there exists a constant $C_m=C_m(C_1, C_2,  P, \omega_g, h_0, T)>0$ such that
\beq
\sup_{t\in[T/2,T)}\left\|h_t\right\|_{C^{m+2}(M,\omega_g, h_0)}\leq C_m.
\label{finite extension spatial estimates}
\eeq
In particular, there exists $\tilde C_m>0$ such that
\beq \|S^{h_t}\|_{C^{m}(M,\omega_g, h_0)}\leq \tilde C_m.\eeq
By the evolution equation $\partial_t h_t = -S^{h_t} + P$, one has
\beq
\sup_{t\in[T/2,T)}\left\|\frac{\p h_t}{\p t}\right\|_{C^m(M,\omega_g, h_0)}\leq C_m'.
\eeq
Moreover, for $s<t$ and $s,t\in[T/2,T)$,
\beq
\left\|h_t-h_s\right\|_{C^m(M,\omega_g, h_0)}
\leq\int_s^t\left\|\frac{\p h_\sigma}{\p\sigma}\right\|_{C^m(M,\omega_g, h_0)}d\sigma
\leq C_m'|t-s|.
\eeq
Since $T<\infty$, this implies that $\{h_t\}_{t\in[T/2,T)}$ is a Cauchy sequence in the $C^m$-norm. In particular,  $h_t\to h_T$ in $C^\infty$ for some smooth Hermitian tensor $h_T$.  By \eqref{finite extension assumptions}, $h_T$ is a Hermitian metric. Therefore, we can extend the flow with initial metric $h_T$.
\eproof

\vskip 2\baselineskip

    \section{A parabolic comparison principle and proof of Theorem~\ref{main2}}
In this section, we establish a parabolic comparison principle that applies to the prescribed Hermitian--Yang--Mills flow \eqref{PHYMF}, the classical Hermitian--Yang--Mills flow, and various other related flows. In particular, we prove Theorem \ref{main2}.

\bdefinition\label{def upper lower null vector condition}
Let $(M,\om_g)$ be a Hermitian manifold, let $E$ be a holomorphic vector bundle on $M$,
and let $\sF:\mathrm{Herm}^+(E)\times \mathrm{Herm}^+(E)\to \mathrm{Herm}(E)$ be a map.
\bd
\item We say that $\sF$ satisfies the \emph{upper null-vector condition} if, for any Hermitian metrics $k$ and $\ell$ such that $k-\ell\geq 0$, one has 
\beq
\sF(k,\ell)(v,\bar v)\geq 0
\label{upper null vector condition}
\eeq
for all nonzero $v$ satisfying $(k-\ell)(v,\bar v)=0$.
\item We say that $\sF$ satisfies the \emph{lower null-vector condition} if, for any Hermitian metrics $k$ and $\ell$ such that $k-\ell\leq 0$, then
\beq
\sF(k,\ell)(v,\bar v)\leq 0,
\label{lower null vector condition}
\eeq
for any nonzero $v$ satisfying $(k-\ell)(v,\bar v)=0$.
\item Fix a smooth Hermitian metric $h_0$ on $E$. We say that $\sF$ is \emph{locally Lipschitz continuous} in the first variable if, for any set $K\subset \mathrm{Herm}^+(E)\times \mathrm{Herm}^+(E)$ that is bounded with respect to $h_0$, there exists a constant $L_{K}=L_K(M, h_0, K,\sF)$ such that for all $(k_1,\ell),(k_2,\ell)\in K$,
\beq
\left|\sF\left(k_1,\ell\right)-\sF\left(k_2,\ell\right)\right|_{h_0}\leq L_K\left|k_1-k_2\right|_{h_0} \quad\text{on }M.
\eeq
Similarly one can define the local Lipschitz continuity in the second variable. \ed
\edefinition
\bexample  Suppose that  $\sF_0(h)= a h+ b P+c^2 \left(\mathrm{tr}_{h_0}h\right) h_0$ where $a, b, c$ are smooth functions on $M$. It is straightforward to verify that
\beq\sF (k,\ell)=\sF_0(k)-\sF_0(\ell)= a (k-\ell)+c^2\left(\mathrm{tr}_{h_0}(k-\ell)\right) h_0\eeq
is locally Lipschitz continuous in both variables and satisfies the lower and upper null-vector conditions.
\eexample

\btheorem\label{strictcomparison} Let $(M,\om_g)$ be a compact Hermitian manifold, and  $E$ be a holomorphic vector bundle on $M$. Suppose that $\underline h(t)$,   $h(t)$ and $\overline h(t)$ are smooth time-dependent Hermitian metrics on $E$ defined for $t\in [0,T)$ with $T\leq \infty$.
\bd \item Let $\sF_>:\mathrm{Herm}^+(E)\times \mathrm{Herm}^+(E)\to \mathrm{Herm}(E)$ be a map satisfying the upper null-vector condition.
If  $\overline h(0)>h(0)$ and
\beq \frac{\p \overline h}{\p t}-\frac{\p h}{\p t}
> S^{h}-S^{\overline h}+\sF_>(\overline h, h),\label{strict upper assumption with F}\eeq
then for any $t\in [0, T)$, one has \beq \overline h(t)>h(t).\eeq
\item Let $\sF_<:\mathrm{Herm}^+(E)\times \mathrm{Herm}^+(E)\to \mathrm{Herm}(E)$ be a map satisfying the lower null-vector condition.
If  $ \underline h(0)< h(0)$ and
\beq \frac{\p \underline h}{\p t}-\frac{\p h}{\p t}
< S^{h}-S^{\underline h}+\sF_{<}(\underline h, h),\label{strict lower assumption with F}\eeq
then for any $t\in [0, T)$, one has \beq  \underline h(t)< h(t).\eeq
\ed

\etheorem

\bproof  We prove only part (1); the proof of part (2) is analogous. Our argument proceeds by contradiction. Define
\beq t_0=\inf\{ t\in [0, T) \ |\ h(t)<\overline h(t) \text{\ fails}\},\eeq
and suppose $t_0<T$. Since $h(0)<\overline h(0)$, by continuity, one deduces that  \beq t_0>0. \eeq  Let \beq Q(t)=\overline h(t)-h(t).\eeq By definition of $t_0$, one concludes that:
\bd\item For each $t\in [0, t_0]$, $Q(t)\geq 0$.
\item At  $t_0$, there exist a point $x_0\in M$ and a nonzero vector $v\in E_{x_0}$ such that
\beq Q(t_0)(v,\bar v)=\overline h_{t_0}(v,\bar v)-h_{t_0}(v,\bar v)=0.\eeq
\item The following inequality holds:
\beq \left.\frac{dQ(t)(v,\bar v)}{dt}\right|_{t=t_0}\leq
0.\label{contra1}\eeq \ed Employing a key technical innovation developed in our previous paper
\cite{XYY24+}, we extend $v$ to a local holomorphic section $V$ of
$E$ over a neighborhood $U$ of $x_0$ satisfying \beq
\left(\nabla^{\overline h_{t_0}} V\right)(x_0)=0. \eeq Consider the
function $f:U\to \R$ given by \beq f=\frac{h_{t_0}\left(V,\bar
V\right)}{\overline h_{t_0}\left(V,\bar V\right)}. \eeq  According
to the setup, we have $f\leq 1$ and $f(x_0)=1$. This implies that
for all $ w\in T^{1,0}_{x_0}M$, \beq \left(\p\bp \log
f\right)(w,\bar w)\leq 0. \eeq On the other hand,  a detailed
calculation shows that \beq \left(\p\bp \log f\right)_{x_0}(w,\bar
w)\geq \frac{R^{\overline h_{t_0}}(w,\bar w,v,\bar
v)}{|v|_{\overline h_{t_0}}^2}-\frac{R^{h_{t_0}}(w,\bar w,v,\bar
v)}{|v|_{h_{t_0}}^2}. \label{key}\eeq 
Since
$|v|^2_{\overline h_{t_0}}=|v|^2_{h_{t_0}}$, one concludes that  for all $w\in
T^{1,0}_{x_0}M$, \beq R^{\overline h_{t_0}}(w,\bar w,v,\bar v)\leq
R^{h_{t_0}}(w,\bar w,v,\bar v). \eeq Taking the trace with respect to
$\omega_g$, one concludes that  at point $(x_0, t_0)$: \beq
S^{\overline h_{t_0}}(v,\bar v)- S^{h_{t_0}} (v,\bar v)\leq 0.\label{upper
curvature inequality}\eeq On the other hand, by the evolution
inequality \eqref{strict upper assumption with F}, \beq \frac{\p
Q}{\p t} =\frac{\p \overline h}{\p t}-\frac{\p h}{\p t}
>S^{h}-S^{\overline h}+\sF_>(\overline h, h).
\eeq
At $(x_0,t_0)$, we have $\overline h_{t_0}-h_{t_0}\geq0$ and $(\overline h_{t_0}-h_{t_0})(v,\bar v)=0$. By the upper null-vector condition,
\beq
\sF_>(\overline h_{t_0}, h_{t_0})(v,\bar v)\geq0.
\eeq
Together with \eqref{upper curvature inequality}, this gives
\beq \left.\frac{dQ(v,\bar v)}{dt}\right|_{t=t_0}>0.\eeq
This is a contradiction to \eqref{contra1}.
\eproof

\noindent By relaxing the hypotheses of Theorem \ref{strictcomparison}, we obtain the following generalization:
\btheorem\label{comparison} Let $(M,\om_g)$ be a compact Hermitian manifold, and  $E$ be a holomorphic vector bundle on $M$. Suppose that $\underline h(t)$,   $h(t)$ and $\overline h(t)$ are smooth time-dependent Hermitian metrics on $E$ defined for $t\in [0,T)$ with $T\leq \infty$.
\bd
\item Suppose that $\sF_\geq$  is locally Lipschitz continuous in the first variable and satisfies the upper null-vector condition.
If  $\overline h(0)\geq h(0)$ and
\beq \frac{\p \overline h}{\p t}-\frac{\p h}{\p t}
\geq S^{h}-S^{\overline h}+\sF_{\geq}(\overline h, h),\label{weak upper assumption with F}\eeq
then for any $t\in [0, T)$, one has \beq \overline h(t)\geq h(t).\eeq
\item
Suppose that  $\sF_\leq$  is locally Lipschitz continuous in the second variable and satisfies the lower null-vector condition.
If  $ \underline h(0)\leq h(0)$ and
\beq \frac{\p \underline h}{\p t}-\frac{\p h}{\p t}
\leq S^{h}-S^{\underline h}+\sF_{\leq}( \underline h, h),\label{weak lower assumption with F}\eeq
then for any $t\in [0, T)$, one has \beq  \underline h(t)\leq h(t).\eeq
\ed

\etheorem

\bproof We prove part $(1)$. Note that the conclusion is local in nature, and we employ a perturbation argument for $\overline h$.  Fix $\tau\in (0,T)$. Since $M\times[0,\tau]$ is compact, there exists a positive constant $C_1$ depending on $\{\overline h(t)\}_{t\in[0,\tau]}$, $\omega_g$ and $M$ such that
\beq
\frac{\p \overline h}{\p t}+S^{\overline h}\geq -C_1\overline h,\quad\text{on $M\times[0,\tau]$.}
\label{upper weak perturbation C1}
\eeq
Fix a Hermitian metric $h_0$ on $E$. Since $\sF_\geq$ is locally Lipschitz continuous in the first variable, there exists a positive constant $C_2$ depending on  $\omega_g$,  $M$, $h_0$ and  the local Lipschitz constant $L_K$ of $\sF_{\geq }$ on
the bounded set
\beq K:=
\left\{\left((1+\mu)\overline h(t), h(t)\right)\in \mathrm{Herm}^+(E)\times \mathrm{Herm}^+(E)\ |\  t\in [0,\tau],\ \mu\in [0,1]\right\},
\eeq  such that for all $0\leq\delta\leq1$,
\beq
\sF_\geq(\overline h, h)-\sF_\geq((1+\delta)\overline h, h)\geq -C_2\delta\overline h,\quad\text{on $M\times[0,\tau]$.}
\label{upper weak perturbation C2}
\eeq
We set
\beq
C_3:=C_1+C_2+1.
\eeq
For any $\eps>0$ satisfying $\eps e^{C_3\tau}<1$,  we define a perturbation metric
\beq
\overline h_\eps(t):=\left(1+\eps e^{C_3t}\right)\overline h(t).
\eeq
It is clear that
\beq
\overline h_\eps(0)>h(0),\qtq{and}
S^{\overline h_\eps}=\left(1+\eps e^{C_3t}\right)S^{\overline h}.
\eeq
In particular,
\beq  S^h-S^{\overline h_\eps}+\sF_\geq(\overline h_\eps, h)= \left(S^h-S^{\overline h}+\sF_\geq(\overline h, h)\right)-\eps e^{C_3t}S^{\overline h}-\sF_\geq(\overline h, h)+\sF_\geq(\overline h_\eps, h),
\eeq
and
\beq \frac{\p \overline h_\eps}{\p t}-\frac{\p h}{\p t}=\frac{\p \overline h}{\p t}-\frac{\p h}{\p t}+\eps e^{C_3t}\left(\frac{\p \overline h}{\p t}+C_3\overline h\right).\eeq
Hence,
\beq \frac{\p \overline h_\eps}{\p t}-\frac{\p h}{\p t}
-\left(S^h-S^{\overline h_\eps}+\sF_\geq(\overline h_\eps, h)\right)= \left(\frac{\p \overline h}{\p t}-\frac{\p h}{\p t}
-\left(S^h-S^{\overline h}+\sF_\geq(\overline h, h)\right)\right)+A\eeq
where
\beq A:=\eps e^{C_3t}\left(\frac{\p \overline h}{\p t}+S^{\overline h}+C_3\overline h\right)
+\sF_\geq(\overline h, h)-\sF_\geq(\overline h_\eps, h).\eeq
By \eqref{upper weak perturbation C1} and \eqref{upper weak perturbation C2}, one has
\beq A\geq\eps e^{C_3t}(C_3-C_1-C_2)\overline h>0.\eeq
Therefore, by \eqref{weak upper assumption with F},
\beq
\frac{\p \overline h_\eps}{\p t}-\frac{\p h}{\p t}
>S^h-S^{\overline h_\eps}+\sF_\geq(\overline h_\eps, h).
\eeq
Applying Theorem~\ref{strictcomparison} to $h$ and $\overline h_\eps$ with $\sF_> = \sF_\geq$ yields
\beq
h(t)<\overline h_\eps(t)=\left(1+\eps e^{C_3t}\right)\overline h(t)
\eeq
for all $t\in[0,\tau]$. Letting $\eps \downarrow 0$ then gives
\beq
h(t)\leq\overline h(t)\label{h(t) leq overline h(t)}
\eeq
for all $t\in[0,\tau]$. Since $\tau \in (0,T)$ is arbitrary, \eqref{h(t) leq overline h(t)} holds for all $t\in [0,T)$. The proof for part $(2)$ is similar.
\eproof

\bremark Theorem~\ref{comparison} provides alternative proofs of the comparison theorems established in \cite{WYY26+} and \cite{FWYY26+}. The framework of Theorem \ref{strictcomparison} and Theorem \ref{comparison} carries  over directly to the setting of Higgs bundles on compact Hermitian manifolds.
\eremark

\noindent \bproof[Proof of Theorem \ref{main2}] Theorem~\ref{main2} follows as a special case of Theorem~\ref{comparison} by taking $\sF_{\geq} = \sF_{\leq} \equiv 0$.
\eproof

\noindent  Theorem \ref{comparison} provides uniform $C^0$-estimates for a general family of flows of the form \begin{equation}
\frac{\partial h}{\partial t} = -\Lambda_{\omega_g}\left(\sqrt{-1}\, R^h\right) + \mathcal{F}_0(h),
\end{equation}
by  defining the comparison function as
\beq \sF(k,\ell)=\sF_0(k)-\sF_0(\ell). \eeq
We  formulate a special case of Theorem~\ref{comparison} directly applicable to the prescribed Hermitian--Yang--Mills flow \eqref{PHYMF}:

\bcorollary\label{cor comparison} Let $(M,\om_g)$ be a compact Hermitian manifold, and  $E$ be a holomorphic vector bundle on $M$. Let $P\in \mathrm{Herm}(E)$ and $\underline h_0, h_0, \overline h_0$ be Hermitian metrics on $E$. Suppose that  $\underline h(t)$, $h(t)$ and $\overline h(t)$ are smooth time-dependent Hermitian metrics on $E$ defined for $t\in [0,T)$ satisfying
\beq \begin{cases} \frac{\p \underline h}{\p t}\leq -S^{\underline h}+P,\\
    \underline h(0)=\underline h_0,\end{cases},   \qquad  \begin{cases} \frac{\p h}{\p t}=-S^h+P,\\
    h(0)=h_0,\end{cases}, \qquad \begin{cases} \frac{\p \overline h}{\p t}\geq -S^{\overline h }+P,\\
    \overline h(0)=\overline h_0.\end{cases}   \eeq
If $\underline h_0\leq h_0\leq \overline h_0$, then for any $t\in [0, T)$,
\beq \underline h (t) \leq h(t)\leq \overline h(t).\eeq
\ecorollary

\noindent As an application of Corollary  \ref{cor comparison}, we establish the long-time existence of the prescribed Hermitian-Yang-Mills flow \eqref{PHYMF}:
\btheorem\label{thm long time existence}
Let  $E$ be a holomorphic vector bundle over a compact Hermitian manifold
$(M,\omega_g)$ and $P\in \mathrm{Herm}(E)$. If there exist two  Hermitian metrics \(\underline h\) and \(\overline h\) such that
\beq
\underline h\leq  \overline h, \qtq{and}
\Lambda_{\omega_g}\left(\sq R^{\underline h}\right)\leq P\leq \Lambda_{\omega_g}\left(\sq R^{\overline
    h}\right),
\eeq
then for any  initial metric $h_0$ on $E$ satisfying  \beq   \underline h\leq h_0\leq \overline h,\eeq  the
prescribed Hermitian-Yang-Mills flow \eqref{PHYMF}
exists for all $t\in[0,\infty)$. Moreover,
\beq \underline h \leq h(t)\leq \overline h, \label{C0}\eeq
holds for all $t\in[0,\infty)$,
\etheorem

\bproof
It is pointed out that the flow \eqref{PHYMF} has a unique smooth solution on a maximal interval $[0,T_{\max})$ and $0<T_{\max}\leq \infty$. On the other hand,  the time-independent metrics $\underline h$ and $\overline h$ are respectively a subsolution and a supersolution of the flow \eqref{PHYMF}, taking themselves as initial metrics.
By Corollary~\ref{cor comparison}, for every $t\in[0,T_{\max})$, one has
\beq
\underline h\leq h(t)\leq \overline h.
\label{global order bound from barriers}
\eeq
 Then there exists a constant $C_1=C_1(h_0,\overline h, \underline h)$ such that $$C_1^{-1}h_0\leq h(t)\leq C_1h_0$$ on $[0,T_{\max})$.  Moreover, there exists a constant $C_2=C_2(C_1, P, h_0)>0$ such that
\beq
\lambda^P_{\min}(x,t)=\inf_{v \in E_x, v\neq 0}  \frac{ P(v, \bar v) }{h(t)(v,\bar v)}\geq -C_2
\label{lambda lower on finite interval}
\eeq
for all $x\in M$ and $t<T_{\max}$. Let $u=|S^{h(t)}-P|_{h(t)}^2$. Then Corollary~\ref{cor tilde S inequality} gives
\beq
\left(\frac{\p}{\p t}-\Delta_\C\right)u\leq 2C_2u.
\label{finite interval tilde S max inequality}
\eeq
 By the parabolic maximum principle, \beq \sup_M u(\cdot,t)\leq e^{2C_2t}\sup_M u(\cdot,0).\eeq   Suppose $T_{\max}<\infty$. Then \beq \sup_M|S^{h(t)}-P|_{h(t)}\eeq  is bounded on $[0,T_{\max})$.  By Theorem \ref{thm finite time extension}, the flow \eqref{PHYMF} extends  over $T_{\max}$ which is absurd. Therefore $T_{\max}=\infty$, and \eqref{C0} holds for all time.   \eproof

\vskip 2\baselineskip

    \section{Proof of Theorem \ref{main3}}
In this section, we prove Theorem \ref{main3}. We begin with the
following lemma which will be essential throughout the proof. Recall
that a matrix-valued function is said to be quasi-positive on
$M$, if it is non-negative definite everywhere, and it is
positive-definite at some point of $M$. \blemma\label{lem main3 DeltaC} Let
$(M,\omega_g)$ be a compact Hermitian manifold and $f$ be a
Lipschitz continuous and quasi-positive  function on $M$.  Then there
exists a positive function $\eta\in C^2(M,\R)$ such that \beq
\left(f-\Delta_\C\right)\eta=1. \label{main3 elliptic weight
concise} \eeq Moreover, if  $f$ is smooth, then $\eta$  is smooth.
\elemma \bproof Fix $p>2\dim_\C M$. It is well-known that
$\Delta_\C:W^{2,p}(M)\to L^p(M)$ is a Fredholm operator of index
zero. Since $W^{2,p}(M)\hookrightarrow L^p(M)$ is a compact
embedding and $f\in L^\infty(M)$, the multiplication map
$W^{2,p}(M)\to L^p(M)$, $u\mapsto fu$, is compact. In particular,
the operator \beq \Delta_\C-f:W^{2,p}(M)\to L^p(M)\eeq is also
Fredholm of index zero. Moreover, if $v \in W^{2,p}(M)$ satisfies
$(\Delta_\C-f)v=0$ in the weak sense, then $v\in
C^{1,\alpha}(M)$ for some $\alpha>0$. Since $f\in C^{0,1}(M)$, by
the Schauder theory for the elliptic equation $\Delta_\C v=fv$, one
concludes that  $v\in C^2(M)$. By the strong maximum principle, we
conclude that $v$ is constant. As $f>0$ somewhere, this constant
must be zero. Hence the kernel of $\Delta_\C-f$ is trivial, and the
index zero property implies that $\Delta_\C-f$ is onto. Suppose that
$\eta\in W^{2,p}(M)$ solves \beq \Delta_\C\eta-f\eta=-1. \eeq
 We also have $\eta\in C^2(M)$. If $x_0$ is a minimum point of $\eta$, then $\Delta_\C\eta(x_0)\geq0$, and so
$
f(x_0)\eta(x_0)=1+\Delta_\C\eta(x_0)\geq1
$.
This implies $\eta>0$ on $M$, and \eqref{main3 elliptic weight concise} follows.
\eproof

\blemma\label{quasipinch} Let  $E$ be a holomorphic vector bundle
over a compact Hermitian manifold $(M,\omega_g)$. If there exists a
smooth Hermitian metric $h$ on $E$ such that \beq
\Lambda_{\omega_g}\left(\sq R^{h}\right) \eeq is quasi-positive,
then for any quasi-positive  $P\in \Gamma\left(M,E^*\ts \bar
E^*\right)$, and any Hermitian metric $h_0$ on $E$,   there exist
two  Hermitian metrics \(\underline h\) and \(\overline h\) on $E$
such that \beq \underline h\leq h_0\leq  \overline h, \qtq{and}
\Lambda_{\omega_g}\left(\sq R^{\underline h}\right)\leq P\leq
\Lambda_{\omega_g}\left(\sq R^{\overline
    h}\right).\eeq

\elemma
\bproof  If  $
\Lambda_{\omega_g}\left(\sqrt{-1}\, R^{h}\right) $ is quasi-positive, we show that there exists a conformal metric $h_E=e^{-\phi} h$ such that $$\Lambda_{\omega_g}\left(\sqrt{-1}\, R^{ h_E}\right)>0.$$
Indeed, it is well-known that there exists some $f\in C^\infty(M,\R)$ such that $\tilde g =e^fg$ and $\omega_{\tilde g}$ is a Gauduchon metric, i.e., $\p\bp\omega_{\tilde g}^{n-1}=0$.  We define
\beq \kappa_{\tilde g, h}(x) = \inf_{0 \neq v \in E_x} \frac{h\left(\Lambda_{\omega_{\tilde g}}\left(\sqrt{-1}\, \Theta^{h}\right) v,\bar v\right)}{h\left(v,\bar v\right)}. \eeq
Recall that $\Theta^h\in \Gamma(M,\Lambda^{1,1}T^*M\ts E^*\ts E)$ is the lifting of $R^h\in \Gamma(M,\Lambda^{1,1}T^*M\ts E^*\ts \bar E^*)$ by $h$.  By using the quasi-positivity and the conformal structure, there exists  some $\lambda_0>0$ such that
\beq \int_M \kappa_{\tilde g, h}(x) \omega_{\tilde g}^n = \lambda_0\int_M \omega_{\tilde g}^n. \eeq
Fix $p>2\dim_\C M$. The elliptic operator
\beq
W^{2,p}(M,\mathbb R)\longrightarrow L^p(M,\mathbb R),
\qquad
u\longmapsto\Lambda_{\omega_{\tilde g}}\left(\sq\p\bp u\right),
\eeq
is Fredholm of index zero. By the maximum principle, its kernel is the one-dimensional space of constants, so its image has codimension one in $L^p(M,\mathbb R)$. The Gauduchon condition implies that this image is contained in the codimension-one subspace
\beq
\left\{f\in L^p(M,\mathbb R):\int_Mf\,\omega_{\tilde g}^n=0\right\}.
\eeq
Since both spaces have codimension one, this inclusion is an equality.
Since $\kappa_{\tilde g,h}$ is Lipschitz and $\lambda_0$ is its $\omega_{\tilde g}^n$-average, there exists a function $\phi\in W^{2,p}(M,\mathbb R)$ such that
\beq \kappa_{\tilde g, h}+\Lambda_{\omega_{\tilde g}}\left(\sq\p\bp\phi\right)= \lambda_0. \eeq
Elliptic regularity gives $\phi\in C^{2,\alpha}(M,\mathbb R)$ for some $0<\alpha<1$.
We define a new metric $h_E= e^{-\phi}h $. By formula \eqref{conformalchange1}, one has
\beq \Lambda_{\omega_{\tilde g}}\left(\sq \Theta^{h_E}\right) - \Lambda_{\omega_{\tilde g}}\left(\sq \Theta^{ h}\right)=  \Lambda_{\omega_{\tilde g}}\left(\sq\p\bp\phi\right) \cdot \mathrm{Id}_E, \label{conformalchange}\eeq
as tensors in $\Gamma(M,E^*\ts E)$.
Moreover, one has
\[ \kappa_{\tilde g, h_E}(x) =  \inf_{0 \neq v \in E_x} \frac{h_E\left(\Lambda_{\omega_{\tilde g}}\left(\sqrt{-1}\, \Theta^{h_E}\right) v,\bar v\right)}{h_E\left(v,\bar v\right)}= \kappa_{\tilde g, h}(x) + \Lambda_{\omega_{\tilde g}}\left(\sq\p\bp\phi\right)=\lambda_0>0. \]
 In particular, the Hermitian-Yang-Mills tensor $\Lambda_{\omega_{\tilde g}}\left(\sqrt{-1}\, \Theta^{h_E}\right)$ is positive definite with respect  to $h^E$. Since $\tilde g =e^f g$, we conclude that $\Lambda_{\omega_g}\left(\sqrt{-1}\, R^{ h_E}\right)>0$. A sufficiently close smooth approximation of $\phi$ in $C^2$ preserves this strict inequality, so we may assume that $h_E=e^{-\phi}h$ is smooth.\\

Since $M$ is compact and $S^{h_E}>0$, there exists a constant $C_1 = C_1(\omega_g, h_0, h_E, P) $ such that
\begin{equation}
h_0 \leq C_1 h_E, \qquad P \leq C_1 S^{h_E}.
\end{equation}
Set $\overline{h} = C_1 h_E$. Then $S^{\overline{h}} = C_1 S^{h_E} \geq P$. Since $P$ is quasi-positive, there exist a non-empty open subset $U \subset M$ and a constant $C_3 = C_3(P, h_E, U) > 0$ such that
\begin{equation}
P \geq C_3 h_E \quad \text{on } U.
\end{equation}
Choose $C_4 > 0$, depending on $\omega_g$ and $h_E$, such that
\begin{equation}
S^{h_E} - C_4 h_E \leq 0 \quad \text{on } M \setminus U.
\end{equation}
Choose $\chi \in C^\infty(M, \R)$ with $\chi \geq 0$, $\chi \not\equiv 0$ and $\supp(\chi)\subset U$. By Lemma~\ref{lem main3 DeltaC}, there exists a smooth function $\psi$ satisfying
\begin{equation}
(\Delta_\C - \chi)\psi = C_4.
\end{equation}
Since $\chi$ is supported in $U$, we have $\Delta_\C \psi = C_4$ on $M \setminus U$, and so
\begin{equation}
S^{h_E} - (\Delta_\C \psi) h_E = S^{h_E} - C_4 h_E \leq 0 \quad \text{on } M \setminus U.
\end{equation}
Consequently, there exists a constant $C_5 = C_5(\omega_g, h_E, \psi)$ such that the following holds on $M$:
\begin{equation}
S^{h_E} - (\Delta_\C \psi) h_E \leq C_5 h_E.
\end{equation}
Furthermore, there exists $C_6 = C_6(h_E, \psi, C_5, h_0, P, U)$ such that
\begin{equation}
C_6 e^\psi h_E \leq h_0 \quad \text{on } M, \qquad
C_6 e^\psi C_5 h_E \leq P \quad \text{on } U.
\end{equation}
Define
\begin{equation}
\underline{h} = C_6 e^\psi h_E.
\end{equation}
The conformal transformation formula \eqref{conformalchange1} for the Hermitian-Yang-Mills tensor yields
\begin{equation}
S^{\underline{h}} = C_6 e^\psi \bigl( S^{h_E} - (\Delta_\C \psi) h_E \bigr).
\end{equation}
By construction, it is straightforward to verify that $S^{\underline{h}} \leq 0 \leq P$ on $M \setminus U$, and $S^{\underline{h}} \leq P$ on $U$. Therefore, $S^{\underline{h}} \leq P$ on $M$, which completes the proof.
\eproof

\bproposition\label{prop uniqueness}
Let $(M,\omega_g)$ be a compact Hermitian manifold, let $E$ be a holomorphic vector bundle on $M$, and  $P\in \mathrm{Herm}(E)$ be quasi-positive. If $h$ and $k$ are two smooth Hermitian metrics on $E$ satisfying
\beq
S^h=S^k=P,
\eeq
then $h=k$.
\eproposition
\bproof
Let $r=\rk(E)$ and set
\beq
\sigma=\tr_h k+\tr_k h.
\eeq
We first compute $\Delta_\C\tr_h k$. Let
\beq
\left(\cE,\tilde h\right)=\left(E^*\ts E,h^*\ts k\right).
\eeq
The identity map $\id_E$ is a holomorphic section of $\cE$, and
\beq
\tr_h k=\left|\id_E\right|_{\tilde h}^2.
\eeq
By the Bochner formula,
\beq
\Delta_\C\tr_h k
= \left|\nabla^{\cE}\id_E\right|_{g,\tilde h}^2
- \LL \Lambda_{\omega_g}\left(\sq \Theta^{\cE}\right)(\id_E),\id_E\RL_{\tilde h}.
\eeq
In a local holomorphic frame $\{e_\alpha\}$ of $E$, using the fact that
\[
\Theta^\cE=\Theta^{E^*,h^*}\ts \id_E+\id_{E^*}\ts \Theta^{E,k},
\]
we obtain
\begin{eqnarray}
\LL \Lambda_{\omega_g}\left(\sq \Theta^{\cE}\right)(\id_E),\id_E\RL_{\tilde h}
&=&g^{i\bar j}\left(-R^h_{i\bar j\alpha\bar\gamma}h^{\beta\bar\gamma}
+R^k_{i\bar j\alpha\bar\gamma}k^{\beta\bar\gamma}\right)
h^{\alpha\bar\delta}k_{\beta\bar\delta}
\nonumber\\
&=&-\LL S^h,k\RL_h+\tr_h S^k.
\end{eqnarray}
Therefore
\beq
\Delta_\C\tr_h k
= \left|\nabla^{\cE}\id_E\right|_{g,\tilde h}^2
+\LL S^h,k\RL_h-\tr_hS^k
\geq \LL S^h,k\RL_h-\tr_hS^k.
\eeq
Similarly,
\beq
\Delta_\C\tr_k h \geq \LL S^k,h\RL_k-\tr_kS^h.
\eeq
Since $S^h=S^k=P$, we obtain
\beq
\Delta_\C\sigma
\geq \LL P,k\RL_h - \tr_hP  + \LL P,h\RL_k - \tr_kP.
\label{uniqueness sigma laplacian}
\eeq
We claim that $\Delta_\C\sigma\geq0$ on $M$. Indeed, fix an arbitrary point
$x\in M$. Choose a basis $\{e_\alpha\}$ of $E_x$ such that at point $x$,
\beq
h_{\alpha\bar\beta}=\delta_{\alpha\beta},\qquad
k_{\alpha\bar\beta}=a_\alpha\delta_{\alpha\beta},\qquad a_\alpha>0.
\eeq
Write $p_\alpha=P_{\alpha\bar\alpha}(x)$. Since $P\geq0$, $p_\alpha\geq0$. At point $x$, a straightforward
computation gives
\beq
\LL P,k\RL_h-\tr_hP
= \sum_{\alpha=1}^r(a_\alpha-1)p_\alpha,\qquad
\LL P,h\RL_k-\tr_kP
= \sum_{\alpha=1}^r(a_\alpha^{-2}-a_\alpha^{-1})p_\alpha.
\eeq
Hence, at point $x$,
\beq
\Delta_\C\sigma
\geq \sum_{\alpha=1}^r \left((a_\alpha-1)+(a_\alpha^{-2}-a_\alpha^{-1})\right)p_\alpha
= \sum_{\alpha=1}^r \frac{(a_\alpha-1)^2(a_\alpha+1)}{a_\alpha^2}p_\alpha
\geq0.
\label{uniqueness sigma subharmonic}
\eeq
Since $x$ is arbitrary, $\Delta_\C\sigma\geq0$ on $M$. By the maximum principle,
$\sigma$ is constant and $\Delta_{\C}\sigma=0$. Let $x\in M$ be a point such that $P(x)$ is positive definite. By \eqref{uniqueness sigma subharmonic},  $a_{\alpha}=1$ for all $\alpha$ since $P_\alpha>0$. This implies $h=k$ at  $x$ and so $\sigma\equiv 2r$. This also implies $h\equiv k$ on $M$.
\eproof

\bproof[Proof of Theorem \ref{main3}]  By Lemma~\ref{quasipinch}, there exist two Hermitian metrics $\underline{h}$ and $\overline{h}$ on $E$ such that
\begin{equation}
\underline{h} \leq h_0\leq \overline{h}, \qquad
\Lambda_{\omega_g}\left(\sqrt{-1}\, R^{\underline{h}}\right) \leq P \leq \Lambda_{\omega_g}\left(\sqrt{-1}\, R^{\overline{h}}\right).
\end{equation}
By Theorem~\ref{thm long time existence}, the solution exists for all $t \geq 0$ and satisfies $\underline{h} \leq h_t \leq \overline{h}$. Therefore, there exists a constant $C_8 = C_8(h_0, \underline{h}, \overline{h}) \geq 1$ such that for all $t \geq 0$,
\begin{equation}
C_8^{-1} h_0 \leq h_t \leq C_8 h_0.
\label{eq:main3_uniform_equivalence_concise}
\end{equation}
We proceed to prove the convergence of the flow. Define
\begin{equation}
u(x,t) := |S^{h_t} - P|_{h_t}^2, \qquad
\lambda_0(x) := \inf_{v \in E_x, \, v \neq 0} \frac{P(v, \bar{v})}{h_0(v, \bar{v})}.
\end{equation}
Then $\lambda_0 \in C^{0,1}(M)$ and $\lambda_0$ is quasi-positive. From $h_t \leq C_8 h_0$, we deduce
\begin{equation}
\lambda^P_{\min}(x,t) := \inf_{v \in E_x, \, v \neq 0} \frac{P(v, \bar{v})}{h_t(v, \bar{v})} \geq C_8^{-1} \lambda_0(x).
\end{equation}
By Corollary~\ref{cor tilde S inequality},
\begin{equation}
\left( \frac{\partial}{\partial t} - \Delta_\C \right) u + 2 C_8^{-1} \lambda_0 u \leq 0.
\label{eq:main3_weighted_inequality_concise}
\end{equation}
Set $f = 2 C_8^{-1} \lambda_0$. Since $f \in C^{0,1}(M)$ is quasi-positive, Lemma~\ref{lem main3 DeltaC} provides a positive function $\eta \in C^2(M)$ satisfying
\begin{equation}
(f - \Delta_\C) \eta = 1.
\end{equation}
Define
\begin{equation}
C_9 := \left( \sup_M \eta \right)^{-1}, \qquad
C_{10} := \sup_M \frac{u(\cdot, 0)}{\eta}.
\end{equation}
Consider the auxiliary function \beq w = u - C_{10} e^{-C_9 t} \eta.\eeq  Clearly $w(\cdot, 0) \leq 0$. Moreover, by \eqref{eq:main3_weighted_inequality_concise},
\begin{equation}
\left( \frac{\partial}{\partial t} - \Delta_\C \right) w + f w \leq C_{10} e^{-C_9 t} (C_9 \eta - 1) \leq 0.
\end{equation}
The parabolic maximum principle implies $w \leq 0$. Consequently, there exists a constant $C_{11}$, depending on $C_{10}$ and $\eta$, such that for all $t \geq 0$,
\begin{equation}
\sup_M |S^{h_t} - P|_{h_t} \leq C_{11} e^{-\frac{C_9}{2} t}.
\label{eq:main3_exponential_decay_concise}
\end{equation}
Since $\partial_t h_t = -(S^{h_t} - P)$, estimates \eqref{eq:main3_uniform_equivalence_concise} and \eqref{eq:main3_exponential_decay_concise} imply
\begin{equation}
\left\| \frac{\partial h_t}{\partial t} \right\|_{C^0(M,\omega_g,h_0)} \leq C_8 C_{11} e^{-\frac{C_9}{2} t}.
\end{equation}
Thus, for $t \geq s \geq 0$,
\begin{equation}
\| h_t - h_s \|_{C^0(M,\omega_g,h_0)} \leq \int_s^t \left\| \frac{\partial h_\sigma}{\partial \sigma} \right\|_{C^0(M,\omega_g,h_0)} d\sigma \leq 2 C_8 C_{11} C_9^{-1} e^{-\frac{C_9}{2} s}.\label{c0exponent}
\end{equation}
Hence $h_t$ converges uniformly to a continuous Hermitian tensor $h_\infty$. Taking the limit in \eqref{eq:main3_uniform_equivalence_concise}, we conclude that $h_\infty$ is a continuous Hermitian metric on $E$.
Finally, utilizing \eqref{eq:main3_uniform_equivalence_concise} and \eqref{eq:main3_exponential_decay_concise}, Theorem~\ref{parabolicbootstrap} implies that for every integer $m \geq 0$, there exists a constant $\tilde{C}_m = \tilde{C}_m(\omega_g, h_0, P, C_8, C_{11})$ such that for all $t \geq 1$,
\begin{equation}
\| h_t \|_{C^m(M,\omega_g,h_0)} \leq \tilde{C}_m.\label{unif}
\end{equation}
The Arzel\`a-Ascoli theorem implies that every sequence $t_i \to \infty$ admits a subsequence converging in $C^\infty$ to a smooth Hermitian metric. Since we have already established $C^0$ convergence, the limit must be $h_\infty$. Thus, the entire flow converges to $h_\infty$ in $C^\infty$. Passing to the limit in \eqref{eq:main3_exponential_decay_concise} yields
\begin{equation}
S^{h_\infty} = P,
\end{equation}
which is equivalent to $\Lambda_{\omega_g}\left(\sqrt{-1} R^{h_\infty}\right) = P$.
Moreover, for all integers $m\geq 1$,  the standard Gagliardo--Nirenberg interpolation inequality implies that 
\beq 
\left\|h_t-h_\infty\right\|_{C^m(M,\omega_g,h_0)}
\leq C_{12}
\left\|h_t-h_\infty\right\|_{C^0(M,\omega_g,h_0)}^{\frac{1}{2}}
\left\|h_t-h_\infty\right\|_{C^{2m}(M,\omega_g,h_0)}^{\frac{1}{2}},
\eeq 
where $C_{12}$ depends on $M, \omega_g,h_0$ and $m$. Moreover, by \eqref{c0exponent} and \eqref{unif},
\beq 
\left\|h_t-h_\infty\right\|_{C^m(M,\omega_g,h_0)}
\leq C_{13}e^{-\frac{C_9}{4}t},
\eeq 
where $C_{13}$ depends on $M, \omega_g,h_0, P$ and $m$. Hence, $\{h_t\}$ converges smoothly  to  $h_\infty$ at an exponential rate.
The uniqueness of $h_\infty$ follows from Proposition \ref{prop uniqueness}. 
 This completes the proof.
\eproof

    \vskip 2\baselineskip

\section{Proof of  Theorem~\ref{main4}}

\noindent\bproof[Proof of Theorem~\ref{main4}]
We first consider $h^-(t)$. By Theorem~\ref{thm long time existence}, it exists for all time and
\beq
\underline h\leq h^-(t)\leq\overline h.
\eeq
The comparison principle first gives $h^-(s)\geq h^-(0)$ for every $s>0$. Since $h^-(s+\bullet)$ is again a solution, another application of the comparison principle gives
\beq
h^-(s+t)\geq h^-(t),
\eeq
and so $\p_t h^-(t)\geq0$. Set $h_t=h^-(t)$ and recall that $\tilde S=S^{h_t}-P$. The flow equation and the preceding monotonicity yield
\beq
\tilde S\cdot h_t^{-1}
=-\frac{\p h_t}{\p t}\cdot h_t^{-1}\leq0,
\label{lower monotone residual sign}
\eeq
To establish convergence using this monotonicity, we define, for two smooth Hermitian metrics $h_0$ and $h_1$ on $E$, the logarithmic determinant functional by
\beq
\cL\left(h_1,h_0\right)
:=\int_M\log\det\left(h_1\cdot h_0^{-1}\right)\omega_g^n.
\label{log determinant functional}
\eeq
Differentiating \eqref{log determinant functional} and using \eqref{lower monotone residual sign}, we obtain
\beq
\frac{d}{dt}\cL\left(h_t,\underline h\right)
=\int_M-\tr_E\left(\tilde S\cdot h_t^{-1}\right)\omega_g^n\geq0.
\label{lower logarithmic determinant variation}
\eeq
Since $h^-(0)=\underline h$ and $h_t\leq\overline h$, for every $T>0$ we have
\beq
0\leq\cL\left(h_T,\underline h\right)
=\int_0^T\int_M-\tr_E\left(\tilde S\cdot h_t^{-1}\right)\omega_g^n dt
\leq\cL\left(\overline h,\underline h\right).
\label{lower logarithmic determinant bound}
\eeq
Set
\beq
\zeta(x,t):=-\tr_E\left(\tilde S\cdot h_t^{-1}\right).
\eeq
Then $\zeta\geq0$, and \eqref{lower logarithmic determinant bound} gives
\beq
\int_0^\infty\int_M \zeta(x,t)\,\omega_g^n dt<\infty.
\label{lower scalar residual integrability}
\eeq
Using \eqref{evolutionS}, $\p_t h_t^{-1}=h_t^{-1}\tilde S h_t^{-1}$, and $\nabla h_t=0$, a straightforward calculation gives
\beq
\frac{\p}{\p t}\left(\tilde S\cdot h_t^{-1}\right)
=g^{i\bar j}D_iD_{\bar j}
\left(\tilde S\cdot h_t^{-1}\right)
-\left(\tilde S\cdot h_t^{-1}\right)
\left(P\cdot h_t^{-1}\right).
\label{monotone residual evolution}
\eeq
Here $D$ is the connection on $E^*\ts E$ induced by the Chern connection of $(E,h_t)$.
Choose $C_1>0$ such that $P+C_1\underline h\geq0$.
Since $h_t\geq\underline h$, we have
\[
P+C_1h_t\geq0,
\qquad
P\cdot h_t^{-1}+C_1\id_E\geq0.
\]
Together with $\tilde S\cdot h_t^{-1}\leq0$, this gives
\[
\tr_E\left(
(\tilde S\cdot h_t^{-1})
(P\cdot h_t^{-1}+C_1\id_E)
\right)\leq0.
\]
Therefore, taking the trace of
\eqref{monotone residual evolution} and using the definition of $\zeta$,
we obtain
\beq
\left(\frac{\p}{\p t}-\Delta_\C\right)\zeta
=\tr_E\left(
(\tilde S\cdot h_t^{-1})(P\cdot h_t^{-1})
\right)
\leq-C_1\tr_E\left(\tilde S\cdot h_t^{-1}\right)
=C_1\zeta.
\label{lower scalar residual differential inequality}
\eeq
Using \eqref{lower scalar residual integrability} and \eqref{lower scalar residual differential inequality}, we claim that
\beq
\lim_{t\to\infty}\sup_M \zeta(\bullet,t)=0.
\label{lower scalar residual uniform decay}
\eeq
Indeed, fix $t\geq1$ and $s\in[t-1,t-1/2]$. By \eqref{lower scalar residual differential inequality},
\beq
\left(\frac{\p}{\p \tau}-\Delta_\C\right)
\left(e^{-C_1(\tau-s)}\zeta(x,\tau)\right)\leq0,
\qquad s\leq\tau\leq t.
\eeq
The smooth uniformly elliptic operator $\Delta_\C$ admits a positive smooth heat kernel with respect to the reference measure $\omega_g^n(y)$; see \cite[Sections~6.2--6.4]{Str08}. Let $H(x,y,\tau)$ denote this heat kernel. Then
\beq
v(x,\tau):=\int_M H(x,y,\tau-s)\zeta(y,s)\,\omega_g^n(y)
\eeq
solves
\beq
\left(\p_\tau-\Delta_\C\right)v(x,\tau)=0,
\qquad v(x,s)=\zeta(x,s).
\eeq

\noindent Define
\beq
w(x,\tau):=e^{-C_1(\tau-s)}\zeta(x,\tau),
\eeq

\noindent The parabolic comparison principle then gives $w(x,\tau)\leq v(x,\tau)$ on $M\times[s,t]$. In particular, at $\tau=t$,
\beq
\zeta(x,t)\leq e^{C_1(t-s)}\int_M H(x,y,t-s)\zeta(y,s)\,\omega_g^n(y).
\label{zeta estimate}
\eeq
Since $H$ is the heat kernel of the fixed operator $\Delta_\C$ on $(M,\omega_g)$ and $t-s\in[1/2,1]$, the smoothness of $H$ and the compactness of $M$ give
\[
C_H:=\sup_{\substack{x,y\in M\\\tau\in[1/2,1]}}
H(x,y,\tau)<\infty.
\]
By construction, $C_H$ depends only on $(M,\omega_g)$. Set $C_2=C_2(C_1,M,\omega_g):=e^{C_1}C_H$.
Integrating \eqref{zeta estimate} in $s$ over $[t-1,t-1/2]$, we obtain
\beq
\sup_M \zeta(\bullet,t)
\leq2C_2\int_{t-1}^{t-1/2}\int_M \zeta(y,s)\,\omega_g^n(y)ds.
\eeq
Combining this estimate with $\zeta\geq0$ and \eqref{lower scalar residual integrability}, we obtain \eqref{lower scalar residual uniform decay}.
Since $-\tilde S\cdot h_t^{-1}$ is non-negative,
\beq
|S^{h_t}-P|_{h_t}
\leq \zeta,
\eeq
and thus $\sup_M|S^{h_t}-P|_{h_t}\to0$.
The order bounds and Theorem~\ref{parabolicbootstrap} now give uniform $C^m$-bounds for $h_t$ on $[1,\infty)$ for every $m$. Monotonicity and the upper barrier give a unique pointwise limit $h^-_\infty$. Every sequence tending to infinity has a smoothly convergent subsequence, and every such cluster limit equals this pointwise limit. Hence the entire flow converges to $h^-_\infty$ in $C^\infty$, and passing to the limit gives
\beq
\Lambda_{\omega_g}\left(\sq R^{h^-_\infty}\right)=P.
\eeq

For $h^+(t)$, set $h_t^+=h^+(t)$ and $\tilde S^+=S^{h_t^+}-P$. The same time-translation argument gives $\p_t h_t^+\leq0$, and therefore
\beq
\tilde S^+\cdot(h_t^+)^{-1}
=-\frac{\p h_t^+}{\p t}\cdot(h_t^+)^{-1}\geq0.
\label{upper monotone residual sign}
\eeq
Differentiating \eqref{log determinant functional} and using \eqref{upper monotone residual sign}, we obtain
\beq
\frac{d}{dt}\cL\left(\overline h,h_t^+\right)
=\int_M\tr_E\left(\tilde S^+\cdot(h_t^+)^{-1}\right)\omega_g^n\geq0.
\label{upper logarithmic determinant variation}
\eeq
Since $h^+(0)=\overline h$ and $h_t^+\geq\underline h$, for every $T>0$ we have
\beq
0\leq\cL\left(\overline h,h_T^+\right)
=\int_0^T\int_M\tr_E\left(\tilde S^+\cdot(h_t^+)^{-1}\right)\omega_g^n dt
\leq\cL\left(\overline h,\underline h\right).
\label{upper logarithmic determinant bound}
\eeq
Set $\zeta^+(x,t):=\tr_E\left(\tilde S^+\cdot(h_t^+)^{-1}\right)$. Similarly, the evolution identity \eqref{monotone residual evolution}, with $(h_t,\tilde S)$ replaced by $(h_t^+,\tilde S^+)$, gives
\beq
\left(\frac{\p}{\p t}-\Delta_\C\right)\zeta^+\leq C_1\zeta^+.
\eeq
Using \eqref{upper logarithmic determinant bound}, the same heat-kernel argument gives $\sup_M \zeta^+(\bullet,t)\to0$ as $t\to\infty$. Since $\tilde S^+\cdot(h_t^+)^{-1}$ is non-negative, this implies $\sup_M|S^{h_t^+}-P|_{h_t^+}\to0$ as $t\to\infty$. Theorem~\ref{parabolicbootstrap} and monotonicity therefore show that $h^+(t)$ converges smoothly to $h^+_\infty$ with
\beq
\Lambda_{\omega_g}\left(\sq R^{h^+_\infty}\right)=P.
\eeq
The comparison principle gives
\beq
\underline h\leq h^-(t)\leq h^+(t)\leq\overline h.
\label{monotone flow order}
\eeq

Passing to the limit in \eqref{monotone flow order} gives $\underline h\leq h^-_\infty\leq h^+_\infty\leq\overline h$. Finally, if $k$ is a stationary solution in this interval, comparison with the constant solution $k(t)=k$ gives
\beq
h^-(t)\leq k\leq h^+(t).
\eeq
Letting $t\to\infty$ proves $h^-_\infty\leq k\leq h^+_\infty$.
\eproof

\vskip 2\baselineskip

\section{Proof of Theorem \ref{main1}}

In this section, we present the proof of Theorem \ref{main1}. A uniform $C^0$-estimate  is obtained in Corollary \ref{cor comparison} as a special case of Theorem \ref{main2}, while long-time existence of the flow is provided by Theorem \ref{thm long time existence}. We shall discuss the long-time convergence in detail. It is ensured by a modified Donaldson functional.
Fix two Hermitian metrics $h_0$ and $h_1$ on $E$.
Let $$h(s):[0,1]\>\mathrm{Herm}^+(E)$$ be a smooth family of Hermitian metrics on $E$  with $h(0)=h_0$ and $h(1)=h_1$. In analogy with the Donaldson functional introduced in \cite{Don85}, 
we define a functional associated with a given Hermitian tensor 
$P \in \mathrm{Herm}(E)$ and two Hermitian metrics $h_0$, $h_1$:
\beq
\cM_P(h_0,h_1)
:=\int_0^1\int_M\tr_E\left(\left(S^h-P\right) h^{-1}\cdot\frac{\p h}{\p s}h^{-1}\right)\om_g^n ds.
\label{path definition MP}
\eeq
Motivated by the approach of \cite{Don85}, we establish the following result.
\blemma\label{Pfunctional}   Let $(M,\omega_g)$ be a compact balanced manifold, i.e.,  $d\omega^{n-1}_g=0$. The quantity $\cM_P(h_0,h_1)$ is independent of the choice of the path $\{h(s)\}_{0\leq s\leq1}$ connecting $h_0$ and $h_1$.  \bd
\item  For every smooth path $\{h_s\}_{0\leq s\leq1}$ starting at $h_0$, one has
\beq
\frac{d}{ds}\cM_P(h_0,h_s)
=\int_M\tr_E\left(
\left((S^{h_s}-P)\cdot h_s^{-1}\right)
\cdot\left(\frac{\p h_s}{\p s}\cdot h_s^{-1}\right)
\right)\om_g^n.
\label{Donaldson first variation}
\eeq

\item Along the prescribed Hermitian-Yang-Mills  flow \eqref{PHYMF}, one has
\beq\frac{d}{dt}\cM_P(h_0,h_t)=
-\int_M\left|S^{h_t}-P\right|_{h_t}^2\om_g^n.
\label{modified Donaldson Lyapunov identity}
\eeq
\ed
\elemma

\bproof Let $h^{(0)}(s)$ and $h^{(1)}(s)$ be two smooth paths in $\mathrm{Herm}^+(E)$ from $h_0$ to $h_1$. For $(s,t)\in[0,1]\times [0,1]$, define
\beq
h(s,t):=(1-t)h^{(0)}(s)+t h^{(1)}(s).
\eeq
Then $h(s,t)\in\mathrm{Herm}^+(E)$ and
\beq
h(s,0)=h^{(0)}(s),\qquad h(s,1)=h^{(1)}(s),\qquad
h(0,t)=h_0,\qquad h(1,t)=h_1.
\eeq
We set
\beq
\Phi_s=\frac{\p h}{\p s}\cdot h^{-1},
\qquad
\Phi_t=\frac{\p h}{\p t}\cdot h^{-1}.
\eeq
We claim that
\beq
\frac{\p}{\p t}\int_M\tr_E(K^h\cdot\Phi_s)\om_g^n
=
\frac{\p}{\p s}\int_M\tr_E(K^h\cdot\Phi_t)\om_g^n.\label{key for path independence}
\eeq
Indeed, one has
\beq
\frac{\partial \Phi_s}{\partial t}
-\frac{\partial \Phi_t}{\partial s}
=\Phi_t\cdot\Phi_s-\Phi_s\cdot\Phi_t
=[\Phi_t,\Phi_s].\label{Phi ts commute}
\eeq
The standard variation formula for the Chern curvature gives
\beq
\frac{\p}{\p t}K^h
=\Lambda_{\om_g}\sq \left(\frac{\p}{\p t}\Theta^h\right)
=\Lambda_{\om_g}\sq\,\bp \p^h\Phi_t,\label{pt Kh}
\eeq
where $\p^h$ and $\bp$ denote the $(1,0)$ and $(0,1)$ parts of the Chern connection on $\mathrm{End}(E)$ induced by $(E,h)$. Indeed, let $\Gamma_{i\alpha}^\beta$ be the Christoffel symbols of the Chern connection of $(E,h)$.  By subtracting the Christoffel symbols of a time-independent reference metric $h_0$, it is clear that $$\Gamma':=\frac{\p}{ \p t}\left(\Gamma_{i\alpha}^\beta dz^i\ts e^\alpha\ts e_\beta\right)$$
is globally defined and
\beq \frac{\p}{\p t} K^h=\Lambda_{\om_g}\sq\,\bp\Gamma'\eeq
On the other hand, it is obvious that
\beq \frac{\p}{ \p t}\left(\Gamma_{i\alpha}^\beta dz^i\ts e^\alpha\ts e_\beta\right)=\p^h\Phi_t.\eeq
Hence, \eqref{pt Kh} follows. By Bochner-Kodaira formula, for every $\phi\in\Gamma\left(M,\mathrm{End}(E)\right)$
\beq
\bp\p^{h} \phi+\p^{h}\bp  \phi =  \phi\cdot \Theta^h-\Theta^h\cdot \phi .
\eeq Hence,
\beq
\frac{\p}{\p s}K^h
=\Lambda_{\om_g}\left(\sq\,\bp  \p^h\Phi_s\right)
=-\Lambda_{\om_g}\left(\sq\,\p^h\bp  \Phi_s\right) +[\Phi_s,K^h].\label{ps Kh}
\eeq
Using \eqref{Phi ts commute}, \eqref{pt Kh}, and \eqref{ps Kh}, we first compute
\begin{eqnarray}
&&\frac{\p}{\p t}\int_M\tr_E\left(K^h\cdot\Phi_s\right)\om_g^n
-\frac{\p}{\p s}\int_M\tr_E\left(K^h\cdot\Phi_t\right)\om_g^n\nonumber\\
&=&\int_M\tr_E\left(\left(\Lambda_{\om_g}\left(\sq\,\bp \p^h\Phi_t\right)\right)\cdot\Phi_s\right)\om_g^n
+\int_M\tr_E\left(\left(\Lambda_{\om_g}\left(\sq\,\p^h\bp \Phi_s\right)\right)\cdot\Phi_t\right)\om_g^n\nonumber\\
&&-\int_M\tr_E\left([\Phi_s,K^h]\cdot\Phi_t\right)\om_g^n
+\int_M\tr_E\left(K^h\cdot[\Phi_t,\Phi_s]\right)\om_g^n.
\label{path independence calculation 1}
\end{eqnarray}
Since $\om_g$ is balanced ($d\om^{n-1}_g=0$), integration by parts yields
\beq
\int_M\tr_E\left(\left(\Lambda_{\om_g}\left(\sq\,\bp \p^h\Phi_t\right)\right)\cdot\Phi_s\right)\om_g^n
=
n\sqrt{-1}\int_M\tr_E\left(\p^h\Phi_t\wedge\bp \Phi_s\right)
\wedge\om_g^{n-1},
\label{first integration by parts}
\eeq
 and
\beq
\int_M\tr_E\left(\left(\Lambda_{\om_g}\left(\sq\,\p^h\bp \Phi_s\right)\right)\cdot\Phi_t\right)\om_g^n
=
n\sqrt{-1}\int_M\tr_E\left(\bp \Phi_s\wedge\p^h\Phi_t\right)
\wedge\om_g^{n-1}.
\label{second integration by parts}
\eeq
Therefore, the sum of the first two terms in \eqref{path independence calculation 1} is
\beq
n\sqrt{-1}\int_M\tr_E\left(
\p^h\Phi_t\wedge\bp \Phi_s
+\bp \Phi_s\wedge\p^h\Phi_t
\right)\wedge\om_g^{n-1}=0.\label{2}
\eeq
On the other hand,
\beq-\tr_E\left([\Phi_s,K^h]\Phi_t\right)
+\tr_E\left(K^h[\Phi_t,\Phi_s]\right)=-\tr_E\left(\Phi_sK^h\Phi_t\right)
+\tr_E\left(K^h\Phi_t\Phi_s\right)
=0.
\label{commutator cancellation}
\eeq
Combining \eqref{path independence calculation 1}, \eqref{2},
and \eqref{commutator cancellation}, we prove \eqref{key for path independence}.
By using \eqref{key for path independence}, we obtain
\begin{eqnarray}
\frac{d}{dt}\int_0^1\int_M\tr_E\left(K^h\cdot\Phi_s\right)\om_g^nds
&=&  \int_0^1 \left(\frac{\p}{\p t}\int_M\tr_E\left(K^h\cdot\Phi_s\right)\om_g^n\right) ds \nonumber\\
&=&  \int_0^1 \left(\frac{\p}{\p s}\int_M\tr_E\left(K^h\cdot\Phi_t\right)\om_g^n\right) ds \nonumber\\
&=& \left(\left.\int_M\tr_E\left(K^h\cdot\Phi_t\right)\om_g^n\right)\right|_{s=0}^{s=1}=0,\label{Donaldson functional part}
\end{eqnarray}
where the last step holds since $\Phi_t\equiv 0$ if $s=0$ or $s=1$.
 On the other hand,
\begin{eqnarray}
&&\int_0^1\int_M\tr_E\left(
\left(-P\cdot h^{-1}\right)\cdot
\left(\frac{\p h}{\p s}\cdot h^{-1}\right)\right)\om_g^n ds= \int_0^1 \left(\frac{\p}{\p s}\int_M\tr_E\left(P\cdot h^{-1}\right)\om_g^n\right) ds\nonumber\\
&=& \int_M\tr_E\left(P\cdot h_1^{-1}\right)\om_g^n-\int_M\tr_E\left(P\cdot h_0^{-1}\right)\om_g^n.
\label{P functional part}
\end{eqnarray}
Therefore $\cM_P(h_0,h_1)$ is independent of the choice of smooth path connecting $h_0$ to $h_1$.
Since \eqref{path definition MP} is invariant under scaling the path,
\beq
\cM_P (h_0,h_s)
:=\int_0^s\int_M\tr_E\left(
\left((S^{h_\tau}-P)\cdot h_\tau^{-1}\right)\cdot
\left(\frac{\p h_\tau}{\p \tau}\cdot h_\tau^{-1}\right)
\right)\om_g^n d\tau.
\eeq
Differentiating with respect to $s$ yields \eqref{Donaldson first variation}.
In particular, along the  flow \eqref{PHYMF}, we have
\beq
\frac{d}{dt}\cM_P(h_0,h_t)
= \int_M\tr_E\left(
\left((S^{h_t}-P)\cdot h_t^{-1}\right)\cdot
\left(\frac{\p h_t}{\p t}\cdot h_t^{-1}\right)
\right)\om_g^n
=-\int_M\left|S^{h_t}-P\right|_{h_t}^2\om_g^n.\eeq
This proves identity \eqref{modified Donaldson Lyapunov identity}.\eproof

\btheorem\label{thm modified Donaldson first variation} Let $(M,\omega_g)$ be a compact Hermitian manifold with $d^*\omega_g=0$. If there is a  constant $C\geq 1$ such that  $C^{-1}h_0\leq h_1\leq Ch_0$, then there exists a constant $C'=C'(\om_g,h_0,P,C)$ such that
\beq
\cM_P(h_0,h_1)\geq-C'.
\eeq
\etheorem

\bproof  We consider the linear path
\beq
h_s=s h_1+(1-s)h_0, \qquad 0\leq s\leq 1,
\eeq
connecting $h_0$ to $h_1$. We also use notations
\beq
H:=h_1\cdot h_0^{-1},\qquad H_s:=sH+(1-s)\id_E,\qquad \Psi_s:=\left(H-\id_E\right)\cdot H_s^{-1}.
\eeq
Then for all $s\in [0,1]$,
\beq
h_s=H_s\cdot h_0,\qquad
\frac{\partial h_s}{\partial s}\cdot h_s^{-1}=\Psi_s.
\eeq
Then the transformation formula \eqref{conformalchange1} gives
\beq
K^{h_s}
= \Lambda_{\om}\left( \sq \Theta^{h_s}\right) = K^{h_0}+\Lambda_{\om_g}\left(\sq \,\bp \left(s(\p^{h_0}H)\cdot H_s^{-1}\right)\right).
\eeq
If we set \beq B_s=s(\p^{h_0}H)\cdot H_s^{-1},\eeq then
\begin{eqnarray}
\int_0^1\int_M \tr_E\left(K^{h_s}\cdot\Psi_s\right)\om_g^n ds
&=&\int_M\int_0^1\tr_E\left(K^{h_0}\cdot\Psi_s\right)\omega_g^n ds\nonumber\\
&&+\int_0^1\int_M\tr_E\left(\left(\Lambda_{\om_g}\left(\sq \,\bp  B_s\right)\right)\cdot \Psi_s\right)\om_g^n ds.
\label{line path calculation}
\end{eqnarray}
For each $\lambda>0$, one has
$
\int_0^1\frac{\lambda-1}{(1-s)+s\lambda}\,ds=\log\lambda$. In particular,
\beq
\int_0^1\Psi_s\,ds=\log H.
\label{linear psi integrates to log}
\eeq
This implies
\beq
\int_M\int_0^1\tr_E\left(K^{h_0}\cdot\Psi_s\right)\omega_g^nds
=\int_M \tr_E\left(K^{h_0}\cdot \log H\right)\omega_g^n.\label{529}
\eeq
Moreover, since $\omega_g$ is balanced, integration by parts yields
\beq
\int_M\tr_E\left(\left(\Lambda_{\om_g}\left(\sq\,\bp  B_s\right)\right)\cdot\Psi_s\right)\om_g^n
=\int_M\tr_E\left(\Lambda_{\omega_g}\left(
\sqrt{-1}B_s\wedge\bp \Psi_s\right)\right)\om_g^n.
\label{linear integration by parts}
\eeq
We show the integrand of the right-hand side is nonnegative. Indeed, fix a point $x$ of $M$. Choose holomorphic coordinates $\{z^i\}$ around $x$ such that $g_{i\bar j}(x)=\delta_{ij}$, and choose an $h_0$-unitary basis $\{e^\alpha\}$ of $E_x$ diagonalizing $H$. With this basis, write
\beq
H=\mathrm{diag}(\lambda_1,\ldots,\lambda_r),
\qquad
H_s=\mathrm{diag}(\nu_1(s),\ldots,\nu_r(s)),
\qquad
\nu_\alpha(s)=s\lambda_\alpha+(1-s).
\eeq
Set $\theta=(\p^{h_0}H)(x)$. Since $B_s=s(\p^{h_0}H)H_s^{-1}$, at point $x$ we have
\beq
(B_s)_{i\alpha}^\beta=\frac{s}{\nu_\beta(s)}\theta_{i\alpha}^\beta.
\eeq
Since $H$ is self-adjoint with respect to $h_0$, the adjoint of $\p^{h_0}H$ with respect to $h_0$ is $\bp H$.
Therefore, at point $x$ we have
\beq
 \bp_{ i} H_\alpha^\beta= \bar{\theta_{i\beta}^\alpha}.
\eeq
A straightforward calculation shows that
\beq
\bp\Psi_s=\left(\bp H\right)H_s^{-1}-s\left(H-\id_E\right)H_s^{-1}\left(\bp H\right)H_s^{-1}.
\eeq
Therefore, at point $x$,
\beq
\p_{\bar i}(\Psi_s)_\beta^\alpha=
\left(\frac{1}{\nu_\alpha(s)}-s\frac{\lambda_\beta-1}{\nu_\alpha(s)\nu_\beta(s)}\right)\p_{\bar i} H_\beta^\alpha
=\frac{1}{\nu_\alpha(s)\nu_\beta(s)}\bar{\theta_{i\alpha}^\beta}.
\eeq
Consequently,
\beq
\tr_E\left(\Lambda_{\omega_g}\left(
\sqrt{-1}B_s\wedge\bp \Psi_s\right)\right)
=\sum_{i,\alpha,\beta}\frac{s}{\nu_\alpha(s)\nu_\beta(s)^2}|\theta_{i\alpha}^\beta|^2\geq0.
\label{linear positive term lower bound}
\eeq
It follows from \eqref{line path calculation}, \eqref{529}, \eqref{linear integration by parts} and \eqref{linear positive term lower bound} that
\beq
\int_0^1\int_M \tr_E\left(K^{h_s}\cdot\Psi_s\right)\om_g^nds\geq
\int_M \tr_E\left(K^{h_0}\cdot \log H\right)\omega_g^n.
\eeq
Since $C^{-1}h_0\leq h_1\leq Ch_0$,  there exists a constant $C_1=C_1(\om_g,h_0,P, C)$ such that
\beq
\int_0^1\int_M \tr_E\left(K^{h_s}\cdot\Psi_s\right)\om_g^nds \geq \int_M\tr_E\left(K^{h_0}\cdot \log H\right)\omega_g^n  \geq  -C_1.
\eeq
On the other hand, there exists
$C_2 = C_2(\omega_g, h_0, P, C)$ satisfying
\begin{eqnarray}
&&\int_0^1\int_M\tr_E\left(-\left(P \cdot h_s^{-1}\right)\cdot \left(\frac{\p h_s}{\p s}\cdot h_s^{-1}\right)\right)\om_g^n ds\nonumber\\
&=&  \int_M\tr_E\left(P \cdot h_1^{-1}\right)\om_g^n-\int_M\tr_E\left(P\cdot  h_0^{-1}\right)\om_g^n
\geq -C_2.
\end{eqnarray}
Hence, we obtain the desired estimate.
\eproof

\bproof[Proof of Theorem \ref{main1}]
By Theorem~\ref{thm long time existence}, the flow \eqref{PHYMF} exists for all time and satisfies
\beq
\underline h\leq h_t\leq \overline h
\label{order bound final theorem}
\eeq
for all $t\in[0, \infty)$.
Let
\beq
u(x,t)=|S^{h_t}-P|_{h_t}^2.
\eeq
By Lemma \ref{Pfunctional} and Theorem~\ref{thm modified Donaldson first variation}, there exists a uniform $C_0=C_0(P, \underline h, \overline h, h_0, \omega_g)>0$ such that for every $T>0$,
\beq
\int_0^T\int_M u(x,t)\,\om_g^n dt
=\cM_P(h_0,h_0)-\cM_P(h_0,h_T)
\leq C_0.
\eeq
Since $u\geq0$, it follows that
\beq
\int_0^\infty\int_M u(x,t)\,\om_g^n dt<\infty.
\label{general scalar residual integrability}
\eeq
Moreover,  by \eqref{order bound final theorem}, there exists a constant $C_1=C_1(P, \overline h, \underline h)>0$ such that \beq \lambda^P_{\min}(x,t)=\inf_{v \in E_x, v\neq 0}  \frac{ P(v, \bar v) }{h(t)(v,\bar v)}\geq-C_1.\eeq  By Corollary~\ref{cor tilde S inequality},
\beq
\left(\frac{\p}{\p t}-\Delta_\C\right)u\leq 2C_1u.
\label{subsolution final theorem}
\eeq
By using  \eqref{general scalar residual integrability} and \eqref{subsolution final theorem} and a similar heat-kernel argument as in the proof of Theorem~\ref{main4},  we obtain
\beq
\lim_{t\to\infty}\sup_M|S^{h_t}-P|^2_{h_t}
=\lim_{t\to\infty}\sup_M u(\bullet,t)=0.
\label{C0 convergence tilde S}
\eeq
By Theorem \ref{parabolicbootstrap} and the Arzel\`a-Ascoli theorem, for every sequence $\{t_i\}$ with $t_i\to\infty$, there exists a subsequence (still denoted by $\{t_i\}$) such that $\{h_{t_i}\}\to h_\infty$ in $C^\infty$, where $h_\infty$ is a  smooth Hermitian metric on $E$. Passing to the limit in \eqref{C0 convergence tilde S} yields $
S^{h_\infty}=P$.
\eproof

\vskip 2\baselineskip
\section{Further generalizations and examples}\label{more}

In this section, we generalize Theorem \ref{main1} and present several examples related to the prescribed Hermitian-Yang-Mills flow.
In the following result, we show that if the condition $S^{\underline h}\leq P\leq S^{\overline h}$ in Theorem \ref{main1} is replaced by the strict inequality $S^{\underline h}<P\leq S^{\overline h}$, then the flow  \eqref{PHYMF} converges globally.

\btheorem\label{prop strict lower whole convergence}
Let $E$ be a holomorphic vector bundle over a compact K\"ahler (or balanced) manifold $(M,\om_g)$, and let $P\in \mathrm{Herm}(E)$. Suppose that there exist two Hermitian metrics $\underline h$ and $\overline h$ on $E$ such that
\beq
\underline h\leq \overline h,
\qquad
S^{\underline h}<P\leq S^{\overline h}.
\label{strict lower proposition assumption}
\eeq
Then for any initial metric $h_0$ satisfying
\beq
\underline h\leq h_0\leq \overline h,
\eeq
the prescribed Hermitian-Yang-Mills flow \eqref{PHYMF} converges in $C^\infty$ to a smooth Hermitian metric $h_\infty$ satisfying
\beq
S^{h_\infty}=P.
\eeq
\etheorem

\noindent One of the key ingredients is the following convexity property of the Hermitian-Yang-Mills curvature tensors.
\blemma \label{convex}
Let $E$ be a holomorphic vector bundle over a Hermitian manifold $(M,\om_g)$.
If $h_0$ and $h_1$ are two Hermitian metrics on $E$, then the Hermitian-Yang-Mills tensors satisfy
\beq
S^{h_0+h_1}\le  S^{h_0}+ S^{h_1},
\eeq
as Hermitian tensors in $\Gamma(M,E^*\ts \bar E^*)$.
\elemma
\bproof Consider the Hermitian holomorphic vector bundle
\begin{equation}
\bigl(\mathcal{E},\tilde{h}\bigr)
:=\bigl(E\oplus E,\,h_{0}\oplus h_{1}\bigr).
\end{equation}
Define the holomorphic bundle map
\begin{equation}
\iota\colon E\longrightarrow \mathcal{E},
\qquad
\iota(v)=(v,v).
\end{equation}
Its image $F:=\iota(E)\subset\mathcal{E}$ is a holomorphic subbundle.
The metric induced by $\tilde{h}$ on $F$ coincides via $\iota$ with
$h_{0}+h_{1}$.  Indeed, for every $x\in M$ and $u,v\in E_{x}$,
\begin{equation}
\tilde{h}\left(\iota(u),\overline{\iota(v)}\right)
= h_{0}\left(u,\bar v\right)+h_{1}\left(u,\bar v\right)
= (h_{0}+h_{1})\left(u,\bar v\right).
\end{equation}
Consequently, for $\sigma=\iota(v)$ we have
\begin{equation}
S^{F}(\sigma,\bar{\sigma})
= S^{\,h_{0}+h_{1}}(v,\bar{v}),
\label{eq:SF-equals-Sh0h1}
\end{equation}
where $S^{F}$ denotes the  Hermitian--Yang--Mills tensor of
$(F,\tilde{h}|_{F})$.
Let $B$ denote the second fundamental form of $F\subset\mathcal{E}$.
By the Chern curvature formula for holomorphic subbundles,
\begin{equation}
S^{F}(\sigma,\bar{\sigma})
= S^{\tilde{h}}(\sigma,\bar{\sigma})
- |B\sigma|^{2}_{g,\tilde{h}}
\le S^{\tilde{h}}(\sigma,\bar{\sigma}).
\label{eq:subbundle-HYM}
\end{equation}
With respect to a product holomorphic frame of $E\oplus E$,
the curvature of $\tilde{h}=h_{0}\oplus h_{1}$ is block diagonal.
Hence, for $\sigma=\iota(v)$,
\begin{equation}
S^{\tilde{h}}(\sigma,\bar{\sigma})
= S^{h_{0}}(v,\bar{v})+S^{h_{1}}(v,\bar{v}).
\label{eq:ambient-HYM}
\end{equation}
Combining \eqref{eq:SF-equals-Sh0h1},
\eqref{eq:subbundle-HYM}, and \eqref{eq:ambient-HYM} yields
\begin{equation}
S^{\,h_{0}+h_{1}}(v,\bar{v})
\;\le\;
S^{h_{0}}(v,\bar{v})+S^{h_{1}}(v,\bar{v})
\end{equation}
for every $x\in M$ and every $v\in E_x$.  This completes the proof.\eproof

\bproof[Proof of Theorem \ref{prop strict lower whole convergence}]
Suppose that $ h_{t_j} \> h_\infty $ is a subsequential \(C^\infty\)-limit
guaranteed  by Theorem~\ref{main1}.
From the previous arguments, we already know that
\begin{equation}
S^{h_\infty} = P,
\qquad
\underline{h}\leq h_\infty \leq  \overline{h}.
\label{subsequential limit order strict lower proof}
\end{equation}
We now prove that the entire flow $\{h_t\}$ converges to $ h_\infty $. By compactness, there exists a constant
$C_1 = C_1(P,\omega_g,\underline{h}) > 0 $ such that
\begin{equation}
P - S^{\underline{h}} \geq C_1\underline{h}.
\label{strict lower S gap}
\end{equation}
Similarly, there exists a constant
$C_2 = C_2(\omega_g,\underline{h}) > 0$  such that
\begin{equation}
S^{\underline{h}} \geq -C_2\underline{h}.
\label{lower bound S underline}
\end{equation}
We set
\begin{equation}
\rho=\min\{C_1/C_2,1/2\},
\label{rho choice lower barrier}
\end{equation}
and define a Hermitian metric
\begin{equation}
\gamma :=(1 - \rho)\,\underline{h}.
\end{equation}
Then, by \eqref{subsequential limit order strict lower proof},
\beq 
h_\infty - \gamma
= (h_\infty - \underline{h}) + \rho\underline{h}
> 0.
\eeq 
Moreover, using \eqref{strict lower S gap},
\eqref{lower bound S underline} and \eqref{rho choice lower barrier},
we obtain
\beq 
P - S^{\gamma}
= P - S^{\underline{h}} + \rho S^{\underline{h}}
\geq C_1\underline{h} - \rho C_2\underline{h}
\geq 0.
\eeq 
Thus
\begin{equation}
\gamma < h_\infty,
\qquad
S^{\gamma} \leq P.
\label{gamma lower S barrier}
\end{equation}
For any \( \varepsilon \in (0,1) \), define
\begin{align}
l_\varepsilon
&:= h_\infty - \varepsilon\,(h_\infty - \gamma),
\label{l_eps def} \\[4pt]
k_\varepsilon
&:= h_\infty + \frac{\varepsilon}{1 - \varepsilon}\,
(h_\infty - \gamma).
\label{k_eps def}
\end{align}
Since $ \gamma < h_\infty$, it follows immediately that
\begin{equation}
l_\varepsilon <h_\infty < k_\varepsilon.
\label{strict sandwich barriers}
\end{equation}
By Lemma \ref{convex}, \eqref{subsequential limit order strict lower proof} and
\eqref{gamma lower S barrier},
\begin{equation}
S^{l_\varepsilon}
\leq  (1 - \varepsilon) S^{h_\infty}
+ \varepsilon S^{\gamma}
\leq  P.
\label{lower barrier inequality}
\end{equation}
On the other hand, from the definition of \( k_\varepsilon \),
\[
h_\infty = (1 - \varepsilon)\,k_\varepsilon + \varepsilon\,\gamma.
\]
Applying Lemma \ref{convex} again yields
\beq 
P = S^{h_\infty}
\leq (1 - \varepsilon) S^{k_\varepsilon}
+ \varepsilon S^{\gamma}
\leq  (1 - \varepsilon) S^{k_\varepsilon}
+ \varepsilon P,
\eeq 
and so
\begin{equation}
P \;\le\; S^{k_\varepsilon}.
\label{upper barrier inequality}
\end{equation}
Combining \eqref{lower barrier inequality} and
\eqref{upper barrier inequality}, we obtain
\begin{equation}
S^{l_\varepsilon}
\leq  P
\leq  S^{k_\varepsilon}.
\label{two stationary barriers}
\end{equation}Since $h_{t_j} \to h_\infty$  in $ C^\infty$ and the inequalities in
\eqref{strict sandwich barriers} are strict, for all sufficiently large \( j \),
\beq 
l_\eps \leq  h_{t_j} \leq  k_\varepsilon.
\eeq 
Applying Theorem \ref{thm long time existence} for the flow \eqref{PHYMF} with initial time $ t_j $,
together with \eqref{two stationary barriers}, gives
\beq 
l_\varepsilon \leq  h_t \leq  k_\varepsilon,
\qquad t \ge t_j.
\eeq
Letting $ \varepsilon \downarrow 0 $ and using
\beq 
l_\varepsilon \to h_\infty,
\qquad
k_\varepsilon \to h_\infty
\quad\text{in } C^0,
\eeq 
we conclude that
\beq 
h_t \to h_\infty
\quad\text{in } C^0.
\eeq 
Finally, by the uniform higher-order estimates furnished by  Theorem \ref{parabolicbootstrap},
\[
h_t \to h_\infty
\quad\text{in } C^\infty.
\]
Thus the entire flow \eqref{PHYMF} converges to \( h_\infty \).
\eproof

For a holomorphic vector bundle $E$, let
$\mathrm{Herm}^{\geq 0}(E)$ denote the set of  nonnegative Hermitian tensors in $\mathrm{Herm}(E)$.  Theorem \ref{main1} has the following refined version:
\btheorem
Let  $E$ be a holomorphic vector bundle over a compact K\"ahler (or balanced) manifold
$(M,\omega_g)$.  For given $P\in\mathrm{Herm}^{\geq 0}(E)$,  if there exist two  Hermitian metrics \(\underline h\) and \(\overline h\) on $E$ such that
\beq
\underline h\leq  \overline h, \qtq{and}
\Lambda_{\omega_g}\left(\sq R^{\underline h}\right)\leq P\leq \Lambda_{\omega_g}\left(\sq R^{\overline
    h}\right),
\eeq
then for any  initial metric $h_0$ satisfying  $  \underline h\leq h_0\leq \overline h$,  the
prescribed Hermitian-Yang-Mills flow \eqref{PHYMF}
exists for all $t\in [0,\infty)$ and converges smoothly to a Hermitian metric $h_\infty$ on $E$ satisfying
\beq
\Lambda_{\omega_g}\left(\sq R^{h_\infty}\right)=P.
\eeq

\etheorem

\bproof
Let \( h_\infty \) be a subsequential limit obtained above.
For every \( \varepsilon \in (0,1) \), there exists \( t_\varepsilon > 0 \) such that
\begin{equation}
(1 - \varepsilon)\,h_\infty
\leq 
h_{t_\varepsilon}
\leq 
(1 + \varepsilon)\,h_\infty.
\label{epsilon closeness to limit}
\end{equation}
Since \( S^{h_\infty} = P \) and \( P \ge 0 \), constant rescalings satisfy
\begin{equation}
S^{(1 - \varepsilon)h_\infty}
= (1 - \varepsilon)\,P,
\qquad
S^{(1 + \varepsilon)h_\infty}
= (1 + \varepsilon)P.
\label{rescaled limit barriers}
\end{equation}
Applying the parabolic comparison principle
(Theorem~\ref{comparison}) from time \( t_\varepsilon \),
with the stationary barriers in \eqref{rescaled limit barriers}, yields
\begin{equation}
(1 - \varepsilon)\,h_\infty
\leq 
h_t
\leq 
(1 + \varepsilon)\,h_\infty,
\qquad
t \ge t_\varepsilon.
\label{whole convergence C0 squeeze}
\end{equation}
Because the solution to the flow is uniquely determined by its initial data,
the evolution of $ h_t$ is independent of the choice of $t_\varepsilon $.
Letting $\varepsilon \downarrow 0 $ in \eqref{whole convergence C0 squeeze},
we obtain
\beq 
h_t \rightarrow h_\infty
\quad\text{in } C^0.
\eeq 
Together with higher-order estimates established in Theorem \ref{parabolicbootstrap},
this $ C^0$-convergence upgrades to
\beq 
h_t \rightarrow h_\infty
\quad\text{in } C^\infty.
\eeq 
Thus the entire flow \eqref{PHYMF} converges smoothly to $ h_\infty$.
\eproof

\bexample\label{example}  In this example, we show that Theorem~\ref{main1} can be applied to establish
existence of solutions to the prescribed Hermitian-Yang-Mills equation
for parameters $P \notin \operatorname{Herm}^{\ge 0}(E)$.\\

 Let $E$ be a holomorphic vector bundle over a compact K\"ahler (or balanced) manifold $(M,\omega_g)$. Suppose that $E$ admits a Hermitian metric $h_{+}$ such that
 \begin{equation}
 S^{h_{+}} > 0.
 \end{equation}
 Fix a point $x_{0}\in M$, and choose a smooth real-valued function $\psi$ satisfying
 \begin{equation}
 \psi(x_{0})=0,
 \qquad
 \Delta_{\mathbb{C}}\psi(x_{0})>0.
 \end{equation}
 Choose constants $0<\varepsilon<1$ and $A>0$ such that
 \begin{equation}
 S^{h_{+}}(x_{0})
 <\varepsilon A\,\bigl(\Delta_{\mathbb{C}}\psi(x_{0})\bigr)\,h_{+}(x_{0}).
 \label{eq:P-negative}
 \end{equation}
 Define
 \begin{equation}
 h_{1}=e^{A\psi}h_{+},
 \qquad
 P=(1-\varepsilon)S^{h_{+}}+\varepsilon S^{h_{1}}.
 \end{equation}
 Then there exists a Hermitian metric $h$ with $S^{h}=P$, although
 \begin{equation}
 P\notin \operatorname{Herm}^{\ge 0}(E).
 \end{equation}
 Indeed, by \eqref{conformalchange}, one has
 \begin{equation}
 S^{h_{1}}=e^{A\psi}\bigl(S^{h_{+}}-A(\Delta_{\mathbb{C}}\psi)h_{+}\bigr).
 \end{equation}
 Evaluating at $x_{0}$ gives
 \begin{equation}
 P(x_{0})
 =(1-\varepsilon)S^{h_{+}}(x_{0})+\varepsilon S^{h_{1}}(x_{0})
 =S^{h_{+}}(x_{0})
 -\varepsilon A\bigl(\Delta_{\mathbb{C}}\psi(x_{0})\bigr)h_{+}(x_{0})
 <0,
 \end{equation}
 which shows that $P\notin \operatorname{Herm}^{\ge 0}(E)$. Define the lower barrier by
 \begin{equation}
 \underline{h}=(1-\varepsilon)h_{+}+\varepsilon h_{1}.
 \end{equation}
 By Lemma~\ref{convex},
 \begin{equation}
 S^{\underline{h}}
 \le (1-\varepsilon)S^{h_{+}}+\varepsilon S^{h_{1}}
 = P.
 \end{equation}
 For the upper barrier, choose a constant $C=C(\underline{h},\omega_g,h_{+},P)>0$ sufficiently large so that
 \begin{equation}
 \underline{h}\le C h_{+},
 \qquad
 P\le C S^{h_{+}}.
 \end{equation}
 Setting $\overline{h}=C h_{+}$, we obtain
 \begin{equation}
 \underline{h}\le \overline{h},
 \qquad
 S^{\underline{h}}\le P\le S^{\overline{h}}.
 \end{equation}
 Finally, Theorem~\ref{main1} implies the existence of a Hermitian metric $h$ on $E$ satisfying $S^{h}=P$.
\eexample

\bexample Let $(E,h_0)$ be a Hermitian vector bundle over a compact Hermitian manifold $(M,\omega_g)$. Suppose that
$$\Lambda_{\omega_g}\left(\sq R^{h_0}\right)=0.$$
\bd \item Let $P=h_0$. Then the prescribed Hermitian-Yang-Mills flow \eqref{PHYMF}
has  the solution
$$h(t)=(1+t)h_0$$
for $t\in[0,\infty)$.
\item  Let $P=-h_0$.  Then the prescribed Hermitian-Yang-Mills flow \eqref{PHYMF} has  the solution
$$h(t)=(1-t)h_0$$
for $t\in[0,1)$; it degenerates at $t=1$.
\ed \eexample

\end{document}